\documentclass[%
 aip,
 amsmath,amssymb,
 reprint,%
]{revtex4-1}

\usepackage{graphicx}% Include figure files
\usepackage{dcolumn}% Align table columns on decimal point
\usepackage{bm}% bold math
\usepackage[utf8]{inputenc}
\usepackage[T1]{fontenc}
\usepackage{mathptmx}
\usepackage{etoolbox}

\usepackage{placeins}
\usepackage{rotating}
\usepackage{svg}
\usepackage{amsmath}
\usepackage{float}
\usepackage{booktabs}
\usepackage{graphicx}

\usepackage{soul}
\usepackage{xcolor}
\newcommand{\ra}[1]{\renewcommand{\arraystretch}{#1}} % Command to adjust row spacin

\makeatletter
\def\@email#1#2{%
 \endgroup
 \patchcmd{\titleblock@produce}
  {\frontmatter@RRAPformat}
  {\frontmatter@RRAPformat{\produce@RRAP{*#1\href{mailto:#2}{#2}}}\frontmatter@RRAPformat}
  {}{}
}%
\makeatother
\begin{document}

\preprint{AIP/123-QED}

\title[Neural ODEs for Stiff Systems: Explicit Exponential Integration Methods]{Training Stiff Neural Ordinary Differential Equations with Explicit Exponential Integration Methods}
% Force line breaks with \\
\author{Colby Fronk}
\affiliation{Department of Chemical Engineering; University of California, Santa Barbara; Santa Barbara, CA 93106, USA}
 \altaffiliation{Correspond to colbyfronk@ucsb.edu}%Lines break automatically or can be forced with \\
 
\author{Linda Petzold}%
\affiliation{Department of Mechanical Engineering; University of California, Santa Barbara; Santa Barbara, CA 93106, USA}%
\affiliation{Department of Computer Science; University of California, Santa Barbara; Santa Barbara, CA 93106, USA}
 \altaffiliation{Correspond to petzold@ucsb.edu}%Lines break automatically or can be forced with \\

\date{\today}% It is always \today, today,

\begin{abstract}
  Stiff ordinary differential equations (ODEs) are common in many science and engineering fields, but standard neural ODE approaches struggle to accurately learn these stiff systems, posing a significant barrier to widespread adoption of neural ODEs. In our earlier work, we addressed this challenge by utilizing single-step implicit methods for solving stiff neural ODEs. While effective, these implicit methods are computationally costly and can be complex to implement. This paper expands on our earlier work by exploring explicit exponential integration methods as a more efficient alternative. We evaluate the potential of these explicit methods to handle stiff dynamics in neural ODEs, aiming to enhance their applicability to a broader range of scientific and engineering problems.  We found the integrating factor Euler (IF Euler) method to excel in stability and efficiency.  While implicit schemes failed to train the stiff Van der Pol oscillator, the IF Euler method succeeded, even with large step sizes.  However, IF Euler's first-order accuracy limits its use, leaving the development of higher-order methods for stiff neural ODEs an open research problem. 
\end{abstract}

\maketitle

\begin{quotation}
Stiff systems of ordinary differential equations (ODEs) describe processes where there are widely varying time scales, with faster dynamics that are stable. Such systems are prevalent in fields like chemistry, biology, and physics, where they model complex phenomena such as chemical reactions or biological processes. Neural ordinary differential equations (neural ODEs) are a type of machine learning model that leverages neural networks to approximate the solutions of differential equations directly from data. However, standard neural ODE approaches often struggle with stiff systems due to stability issues, limiting their broader applicability in scientific and engineering research. In previous work, we utilized single-step implicit schemes to address this challenge, enabling neural ODEs to effectively handle stiffness. While this approach was successful, implicit methods have a high computational cost. In this paper, we explore the use of explicit exponential integration methods, as a more efficient alternative for training stiff neural ODEs. These methods maintain stability while reducing computational overhead, potentially expanding the range of complex problems that neural ODEs can tackle in scientific and engineering domains.

\end{quotation}

\section{Introduction}

No discipline can fully grasp the intricacies of complex systems without the rigorous development of mathematical models, which serve as a cornerstone for advancing our understanding of natural and engineered systems.  For instance, in epidemiology, ordinary differential equations (ODEs) model the spread of infectious diseases such as influenza, measles, and COVID-19, while in medicine, they describe the dynamics of CD4 T-cells in HIV patients. Similarly, partial differential equations (PDEs) are fundamental in climate science, enabling the modeling of atmospheric and oceanic processes to forecast and mitigate risks such as severe droughts, flooding, and extreme storms.  Crafting a mathematical model with adequate detail is crucial, as it enables the identification of effective intervention strategies, such as pharmaceutical treatments, to mitigate undesirable outcomes such as the spread of disease.  In engineering, mathematical models facilitate the design of control systems, ensuring reliability and stability in diverse conditions.  However, the traditional model development cycle is time-intensive, involving finding a candidate model to describe processes, using data to fit parameters to the model, analyzing uncertainties in the fitted parameters, and performing additional experiments to refine and validate the model. 

System identification has been revolutionized by Sparse Identification of Nonlinear Dynamics (SINDy) \cite{SINDY, BayesianSINDy, Kaheman2020SINDyPIAR}, a regression technique applied on numerical approximations to the derivative.  SINDy has been successfully applied to ODE model extraction in diverse domains such as fluid dynamics \cite{PDE_SINDy}, plasma physics \cite{plasma_SINDy}, and chemical reaction networks \cite{reaction_networks_SINDy, reactive_SINDy}.  However, its reliance on densely sampled data limits broader adoption \cite{doi:10.1063/5.0130803}.

With unprecedented data from IoT devices \cite{li2015internet, rose2015internet}, automated high-throughput biological research systems \cite{trends_highthroughput_screening, szymanski2011adaptation}, and Earth observation systems \cite{satellite}, data-driven modeling has emerged as a powerful approach to identify system dynamics directly without relying on first principle models.  Recent frameworks such as neural ODEs \cite{NeuralODEPaper, latent_ODEs, bayesianneuralode, stochastic_neural_ode, kidger2020neural, kidger2022neural, morrill2021neural, jia2019neural, chen2020learning, dagstuhl, doi:10.1063/5.0130803, 10.1371/journal.pcbi.1012414} , physics-informed neural networks (PINNs) \cite{owhadi2015bayesian, hiddenphysics, raissi2018numerical, raissi2017physics, osti_1595805, cuomo2022scientific, cai2021physics}, and MeshGraphNets \cite{pfaff2021learningmeshbasedsimulationgraph} offer flexible solutions for modeling dynamical systems across diverse scientific and engineering domains.  Scientists and engineers find it challenging to rely on "black box" models due to their lack of interpretability and the need for a clear understanding of how and why they function.  Symbolic neural networks were introduced to overcome these challenges.  By integrating symbolic mathematical expressions into their architecture, symbolic neural networks facilitate symbolic regression and produce interpretable equations that uphold the reliability and precision demanded by scientists and engineers \cite{PiNetPaper, Integration_of_Neural, Kubal_k_2023, zhang2023deep, 9659860, doi:10.1063/5.0130803, su2022kinetics, Ji_2021, 10.1063/5.0134464, 10.1371/journal.pcbi.1012414}.

Despite their successes, neural ODEs face significant obstacles in handling stiff differential equations, which arise frequently in multiscale phenomena with widely varying timescales.  Explicit methods are inexpensive to implement and use with neural differential equations, but the small number of time steps that stiff system of equations demand for numerical stability makes them too expensive to use with stiff neural differential equations.  Even when training data is derived from non-stiff ODEs, neural ODEs can become stiff during training due to the highly nonlinear dynamics introduced by neural network models. This unintended stiffness can drastically slow down training or even prevent convergence, necessitating that all neural ODE solvers be inherently robust to stiffness.  While some approaches have attempted to mitigate stiffness in neural ODEs by modifying the system dynamics, such as through equation scaling or regularization \cite{stiff_neural_ode, caldana2024neural, dikeman2022stiffness, LINOT2023111838, holt2022neural, baker2022proximal, MALPICAGALASSI2022110875, thummerer4819144eigen, weng2024extending, ghosh2020steer, finlay2020trainneuralodeworld, kelly2020learningdifferentialequationseasy, onken2020discretizeoptimizevsoptimizediscretizetimeseries, onken2021otflowfastaccuratecontinuous, massaroli2020stableneuralflows, massaroli2021dissectingneuralodes, ji2021stiff, guglielmi2024contractivity, pal2021opening, pal2023locallyregularizedneuraldifferential, kumara2023physics, massaroli2020stable}, these methods address the problem only indirectly and do not fully resolve the actual stiffness issue.  Our earlier work (see Ref.~\onlinecite{fronk2024trainingstiffneuralordinary}) demonstrated that single-step implicit methods such as backward Euler, trapezoid method, and Radau3 and Radau5 can effectively train stiff neural ODEs without altering their dynamics through scaling or regularization.  Although implicit methods handle stiffness effectively, solving nonlinear systems at every time step can lead to high computational costs.  Given these limitations, an open question emerges: can we develop more efficient neural ODE methods for handling stiffness that avoid the need for nonlinear solvers?

Explicit exponential integration methods \cite{cox2002exponential, hochbruck2010exponential, iyiola2018exponential, maday1990operator, li2023low, krogstad2005generalized, owren1999runge, crouch1993numerical, celledoni2003commutator, tokman2006efficient, certaine1960solution} present a promising direction as a nontraditional class of numerical techniques that manage stiff systems while retaining the simplicity of explicit schemes. Unlike standard explicit methods, these approaches use matrix exponentiation to stabilize the integration process, potentially enabling larger time steps without sacrificing accuracy. Since they do not require solving nonlinear systems at each step, explicit exponential integrators can be more computationally efficient and are inherently differentiable, making them well-suited for direct backpropagation in neural ODE frameworks. 

In this paper, we propose a new method for training stiff neural ODEs by employing the explicit exponential integrating factor Euler (IF Euler) method.  This work builds on our prior work (see Ref.~\onlinecite{fronk2024performanceevaluationsinglestepexplicit}) that conducted an extensive evaluation of higher-order exponential integration methods, systematically comparing their stability and accuracy with widely used implicit schemes.  That study showed that higher-order exponential methods, while innovative, fail to outperform the first-order explicit IF Euler method without sacrificing stability, making them unsuitable for stiff systems. These results led us to conclude that IF Euler is the only viable option for repeated, low-cost integration in neural ODEs, which is why this paper focuses exclusively on this method.

By investigating explicit exponential integration methods, we establish a foundation for designing efficient and stable ODE solvers specifically suited for handling stiffness in neural ODEs. This approach is pivotal for expanding the use of these solvers to diverse and challenging problems, such as partial differential equations (PDEs) with stiffness due to diffusion effects, mesh-based simulations with MeshGraphNets, physics-informed neural networks (PINNs), and other innovative frameworks that demand both differentiability and computational efficiency in solving differential equations.

\FloatBarrier
\section{Methods}
\FloatBarrier
\subsection{Neural ODEs}

\begin{figure*}
    \centering
    \includegraphics[width=0.9\linewidth]{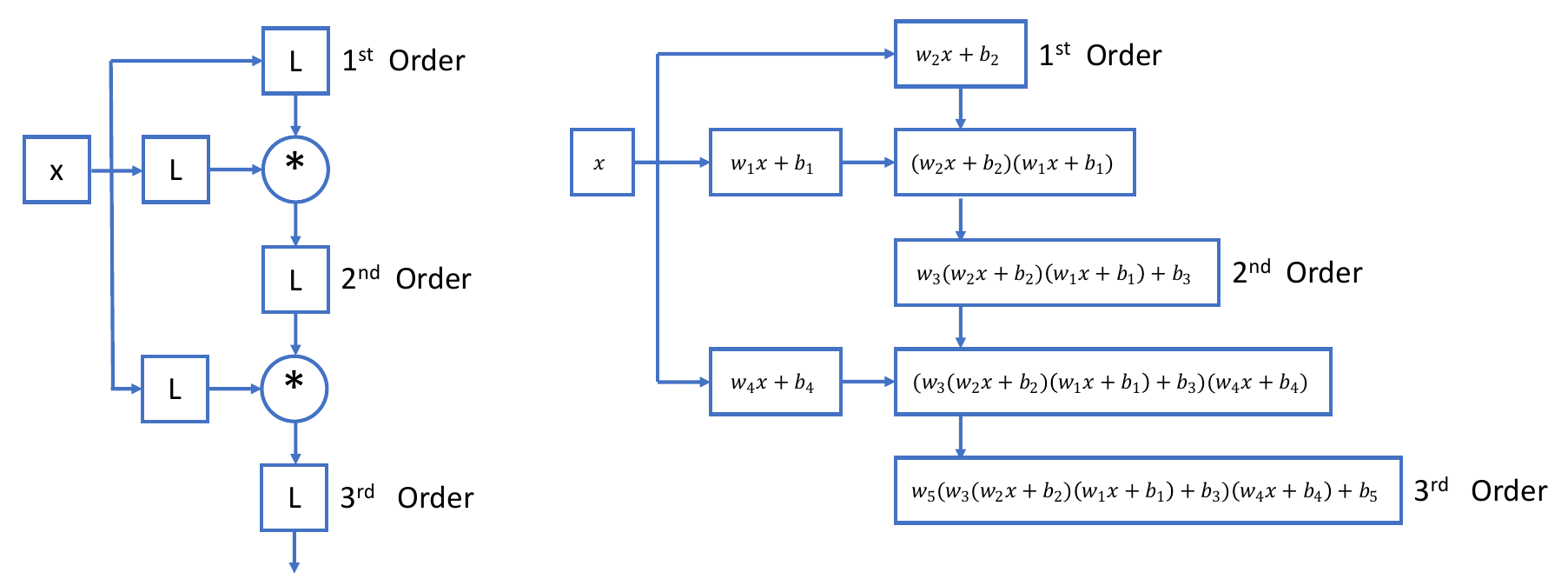}
    \caption{The neural network architecture of $\pi$-net V1 from Ref.~\onlinecite{PiNetPaper} is depicted on the left. On the right, a worked example illustrates how a 1-dimensional input layer with the variable $x$ propagates through the network. Circles marked with the $*$ symbol denote layers where the Hadamard product is applied to the layer's inputs. The boxes labeled $L$ represent standard linear layers without activation functions. This architecture does not utilize typical activation functions like tanh or ReLU, which enhances its interpretability.}
    \label{fig:PiNetArch}
\end{figure*}

Neural Ordinary Differential Equations \cite{NeuralODEPaper} (neural ODEs) are a type of neural network designed to approximate time-series data, $y(t)$, by modeling it as an ODE system. In many scientific fields, the ODE system we aim to approximate takes the following form:

\begin{equation} \frac{dy(t)}{dt} = f\left(t, y(t), \theta \right), \end{equation}

\noindent where $t$ represents time, $y(t)$ is a vector of state variables, $\theta$ denotes the parameter vector, and $f$ is the function defining the ODE model. Finding an accurate form for $f$ can be a complex and time-consuming task.  To address this, neural networks can be used to approximate $f$ by leveraging the universal approximation theorem \cite{Hornik1989MultilayerFN}, allowing us to replace $f$ with a neural network model, $NN$:

\begin{equation} \frac{dy(t)}{dt} = f \approx NN\left(t, y(t), \theta \right). \end{equation}

Like traditional ODE solvers, neural ODEs are solved by integrating from an initial condition using standard differential equation discretization schemes \cite{ascher1998computer, griffiths2010numerical, hairer2008solving} to produce time-series predictions. Most of the time, some parts of the model, denoted as $f_{known}$, are known, but not all mechanisms and terms describing the full model are understood.  Neural ODEs can refine or improve the model by learning only the missing or poorly understood components:

\begin{equation} \label{eqn
} \frac{dy(t)}{dt} = f_{known}\left(t, y(t), \theta \right) + NN\left(t, y(t), \theta \right). \end{equation}

\noindent This method enables the identification and modeling of gaps in existing models, allowing for targeted refinement where traditional modeling is incomplete or uncertain. It combines the strengths of known physics with data-driven discovery to enhance the overall accuracy and completeness of the system representation.

\subsection*{Polynomial Neural ODEs}

Differential equations with polynomial right-hand side functions $f$ commonly appear in fields like chemical kinetics \cite{soustelle2013introduction}, cell signaling \cite{gutkind2000signaling}, gene regulation \cite{peter2020gene}, epidemiology \cite{magal2008structured}, and ecology \cite{mccallum2008population}. When it is known that a system is governed by polynomial relationships, polynomial neural ODEs become particularly effective for solving these inverse problems.

Polynomial neural networks \cite{PiNetPaper, FAN2020383} are a type of symbolic network architecture where the outputs are polynomial functions of the inputs. There are various forms of polynomial neural networks, each with its unique structure and properties. For an extensive overview, readers can consult Grigorios G. Chrysos's work on these architectures. In our research, we found that the $\pi$-net V1, as proposed in Ref.~\onlinecite{PiNetPaper}, was particularly effective (see Fig. \ref{fig:PiNetArch}). This architecture is constructed using Hadamard products \cite{horn1994topics} of linear layers without activation functions:

\begin{equation} L_i(x) = x*w_i+b_i \end{equation}

\noindent where $L_i$ represents linear transformations. This formulation allows for the creation of higher-degree polynomial expressions, and the network’s structure must be predefined according to the desired polynomial degree. Importantly, the $\pi$-net does not require hyperparameter tuning, and it can generate any n-degree polynomial for the given input variables. The dimensions of hidden layers can vary, provided that the Hadamard product operations remain dimensionally compatible. By integrating polynomial neural networks into the neural ODE framework \cite{NeuralODEPaper}, Polynomial Neural Ordinary Differential Equations (PNODEs) \cite{doi:10.1063/5.0130803} offer a novel approach to symbolic modeling. Since the output of a polynomial neural ODE is constructed solely through tensor and Hadamard products without nonlinear activation functions, the final trained network can be directly translated into a symbolic mathematical expression. This is a significant advantage over traditional neural networks and neural ODEs.

Polynomial neural ODEs are particularly beneficial for learning stiff ODEs and evaluating numerical methods due to several key advantages. They do not require data normalization or standardization, making them straightforward to apply across different modeling scenarios. Additionally, their ability to manage outputs over a wide range of scales is crucial for stiff ODEs, which often involve constants spanning multiple orders of magnitude, offering a robust solution for such complex systems.

\begin{figure*}
    \centering
    \includegraphics[width=0.95\linewidth]{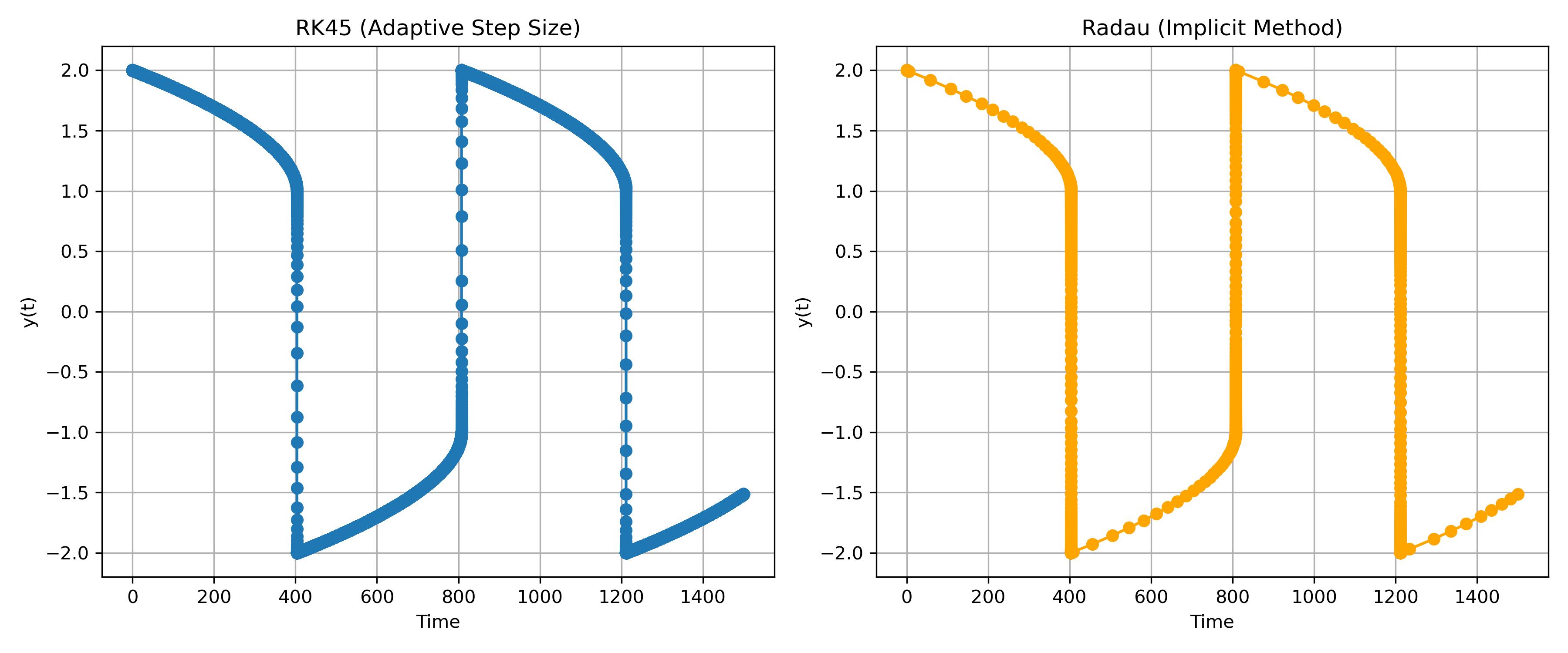}
    \caption{Comparison of the integration of the deterministic stiff van der Pol oscillator with $\mu=1000$ using two different methods: (a) explicit Runge-Kutta-Fehlberg, which is slow with 422,442 time points and 2,956,574 function evaluations, and (b) implicit Radau IIA 5th order, which is faster with only 857 time points and 7,123 function evaluations.}
    \label{fig:stiffness_explained}
\end{figure*}

\subsection{Stiff ODEs}

Stiff ODEs present a challenge where classical explicit integration methods become inefficient due to stability constraints. Stiffness occurs when there is a large disparity in time scales, often indicated by widely varying eigenvalues in the Jacobian matrix. For example, the stiff van der Pol oscillator with $\mu = 1000$ (Figure 2) requires 422,442 data points for stable integration using the explicit Runge-Kutta-Fehlberg method, making the computation slow and costly. In contrast, the implicit Radau IIA 5th order method reduces this to 857 points, significantly speeding up the process. However, the computational cost of implicit methods isn't fully captured by the number of time points, as they involve iterative processes. A better measure is the number of function evaluations: the explicit method required 2,956,574 evaluations, whereas Radau IIA required only 7,123, leading to much faster integration.

This example illustrates the core issue of stiffness in differential equations—classical explicit methods require very small time steps for stability, greatly increasing computational cost. This is especially problematic for neural ODEs, where thousands of integrations are needed during training. Although classical explicit methods are simpler and cheaper per step, their inefficiency in stiff problems makes them impractically slow.

\begin{figure*}
    \centering
    \includegraphics[width=0.9\linewidth]{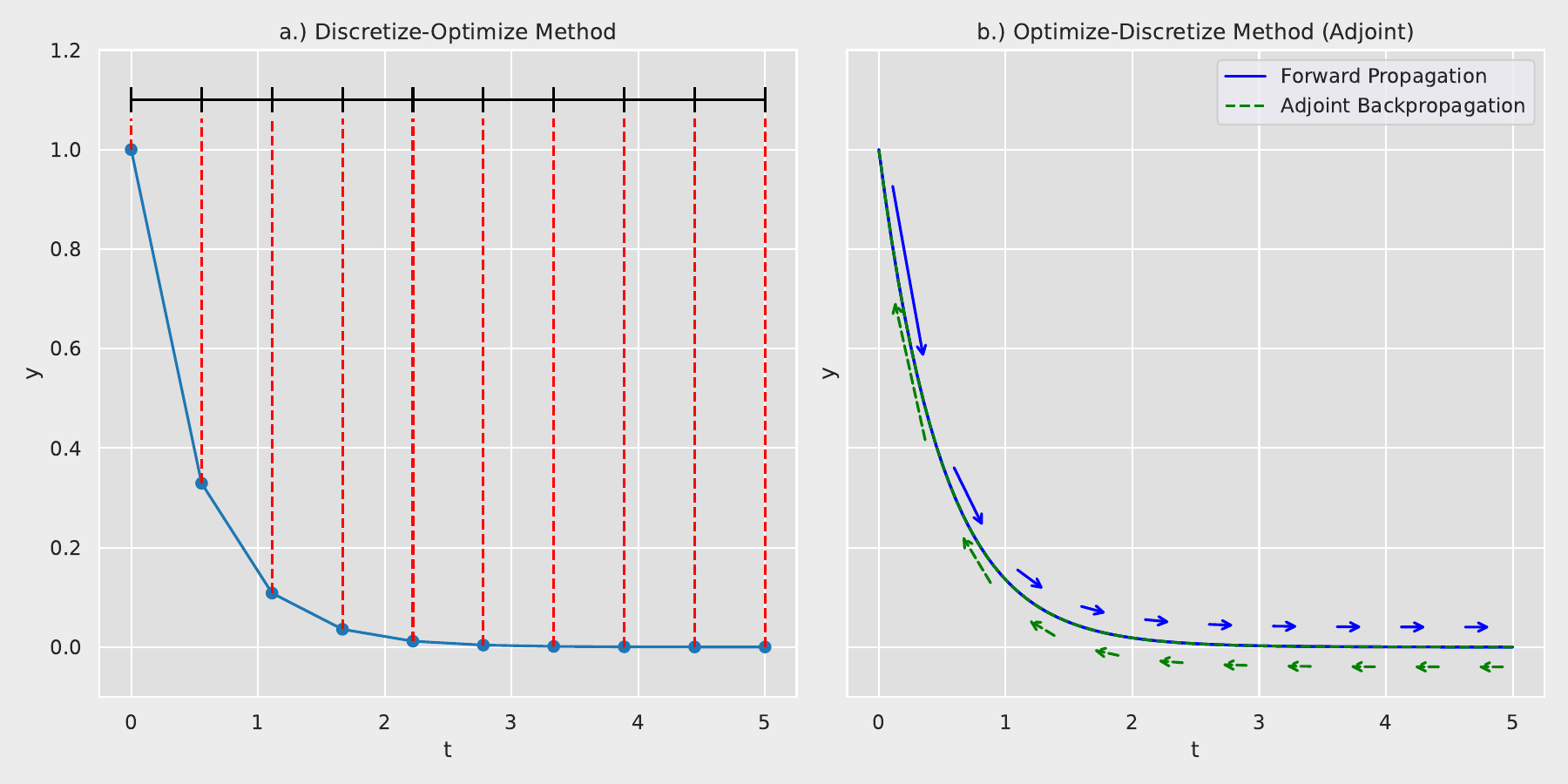}
    \caption{Illustration of (a) Discretize-Optimize and (b) Optimize-Discretize methods. For Discretize-Optimize, black and red lines denote the discretized grid used to perform the optimization. For Optimize-Discretize, blue arrows indicate the forward pass of the neural network, while blue lines depict the backward pass using the adjoint method, illustrating how gradients are computed.}
    \label{fig:Explaining_Disc-Opt_vs_Opt-Disc}
\end{figure*}

\subsection{Discretize-Optimize vs. Optimize-Discretize}

Optimizing neural ODEs generally involves two primary strategies: Discretize-Optimize (Disc-Opt) and Optimize-Discretize (Opt-Disc). These methods differ fundamentally in their approach to solving and optimizing ODE problems, as depicted in Fig. \ref{fig:Explaining_Disc-Opt_vs_Opt-Disc}. Disc-Opt starts by discretizing the ODE, allowing optimization to be performed directly on this discretized form, which is easier to implement, especially with automatic differentiation tools for efficient gradient computation \cite{wright2006numerical}. Conversely, Opt-Disc defines an optimization problem in the continuous domain and computes gradients before discretizing, requiring the numerical solution of the adjoint equation \cite{bliss1919adjoint}, which is computationally more demanding.

Both strategies have their pros and cons, particularly concerning the computational burden of solving the ODE during forward propagation and calculating gradients for optimization. Solving forward propagation accurately requires considerable memory and computational power, which can be prohibitively expensive during training. While using a less accurate solver might speed up the process, it introduces risks for the Opt-Disc method, where inaccuracies in forward and adjoint solutions can lead to poor gradient quality \cite{gholami2019anode}. This issue does not affect Disc-Opt, where gradient accuracy is decoupled from the forward solver's precision, providing more flexibility in optimizing solver accuracy and improving training efficiency. Opt-Disc's reliance on adjoint methods for backward computation of the neural ODE also makes it more prone to numerical instabilities, particularly with stiff ODEs, impacting the robustness and effectiveness of training.

Research by Onken et al. \cite{onken2020discretizeoptimizevsoptimizediscretizetimeseries} supports the benefits of the Disc-Opt approach. It shows that Disc-Opt not only achieves similar or better validation loss outcomes compared to Opt-Disc but also does so with significantly reduced computational costs—yielding an average speedup of 20x. This efficiency arises partly because Opt-Disc methods can produce unreliable gradients if the state and adjoint equations lack sufficient precision \cite{gholami2019anode}. In contrast, Disc-Opt maintains gradient reliability regardless of solver accuracy, allowing for adjustments based on data noise levels, which is advantageous in scientific modeling.

Given these findings, our research focuses on the Discretize-Optimize approach to develop a more efficient differential ODE solver. By building on the advantages of Disc-Opt, particularly for stiff ODEs, we aim to refine training strategies using different integration schemes, ultimately creating solvers that are both robust and computationally efficient.

\subsection{Classical Numerical Methods for Solving ODEs}

Solving ordinary differential equations (ODEs) typically involves two main types of methods: single-step and multistep approaches. Single-step methods, such as Euler and Runge-Kutta, calculate the solution at the next time step using only the current state. This simplicity makes them easier to implement and analyze since they do not require any historical data from previous steps. On the other hand, multistep methods, including Adams-Bashforth and Adams-Moulton, derive the next value based on several preceding time steps. While these methods can offer computational efficiency in some scenarios, they are more complex to manage due to their dependence on past values and their intricate stability requirements. Thus, multistep methods, while potentially more efficient, involve added complexities not present in single-step methods.

Single-step methods are favored for their simplicity and ease of implementation, making them a natural choice for initial studies in neural ODE integration. Because they depend only on the current state, they simplify the process of understanding neural ODE dynamics and enable straightforward backpropagation through the computed ODE solution. Conversely, multistep methods pose specific difficulties in the context of neural ODEs. These methods require accurate data from several previous steps, which can be problematic when training data is noisy or incomplete. In such cases, past points needed for the method must be recalculated, increasing computational effort and complexity. Moreover, backpropagation with multistep methods is more challenging because gradients must be propagated not just through the current step but also through the reconstructed previous steps. Due to these complications, single-step methods should be explored during the early stages of neural ODE research, while more complex multistep methods are better suited for later exploration once a solid understanding of single-step methods has been established.

\subsection*{Implicit Methods are Expensive}

Implicit methods offer a robust approach to handling the difficulties associated with stiff ODEs. In contrast to explicit methods, which determine the next step using only current state values, implicit methods involve solving a system of nonlinear equations at every time step. Although this requires more computational effort, it delivers much higher stability, allowing for larger time steps while maintaining accuracy. This balance between computational expense and stability makes implicit methods ideal for stiff problems; however, they can be less efficient for non-stiff ODEs where stability is not a concern.

Some common implicit methods for solving stiff ODEs are the backward Euler method, the trapezoid method, and Radau IIA methods. The backward Euler method is the simplest approach:  

\begin{equation}
y_{n+1} = y_n + h f(t_{n+1}, y_{n+1}),
\end{equation}

\noindent Backward Euler is A-stable, providing strong damping of oscillations and stability for stiff systems. However, it can overdampen solutions, making them too smooth. The trapezoid method improves upon backward Euler by averaging the derivatives at the current and next points: 

\begin{equation}
y_{n+1} = y_n + \frac{h}{2} \left( f(t_n, y_n) + f(t_{n+1}, y_{n+1}) \right).
\end{equation}

\noindent Trapezoid method is second-order accurate and A-stable but not L-stable, meaning it may allow some stiff components to oscillate rather than decay smoothly. While it offers higher accuracy, the potential for introducing oscillations can be problematic for stiff neural ODEs during training. Radau IIA methods \cite{axelsson1969class, ehle1969pade, hairer1999stiff} have the following form:

\begin{align}
Y_i =& y_n + h \sum_{j=1}^{s} a_{ij} f(t_n + c_j h, Y_j), \quad i = 1, \ldots, s, \\
y_{n+1} =& y_n + h \sum_{j=1}^{s} b_{j} f(t_n + c_j h, Y_j)
\end{align}

\noindent where \(Y_i\) are the stage values, \(a_{ij}\) are the coefficients from the Butcher tableau, and \(c_j\) are the nodes. The Butcher tableau for Radau3 and Radau5 can be found in Ref.~\onlinecite{hairer1999stiff}. The Radau3 and Radau5 methods are both A-stable and L-stable, which ensures rapid decay of transient components without oscillations, even with large time steps. 

For these methods, the solution at the next step, $y_{n+1}$ is found by solving a multivariate nonlinear system of equations using techniques such as Newton's method.  This adds significant computational cost and may require step size adjustments if convergence fails.  Although the Radau IIA methods offer greater accuracy, they require solving for \(Y_i\) a nonlinear solution involving a matrix-valued function, which significantly increases computational cost and introduces additional convergence challenges.

Backpropagation is the process of tracing a computational graph to compute gradients of a loss function with respect to model parameters. For neural ODEs, explicit integration methods like Euler and RK4 compute future values directly from the current state, forming a straightforward computational graph that is easy to backpropagate through. In contrast, implicit methods like the backward Euler or Radau IIA methods solve a nonlinear equation at each step, creating a nested graph structure. This makes backpropagation computationally expensive and susceptible to numerical issues when unrolling these loops. The implicit function theorem provides a workaround by computing gradients directly at the solution, avoiding the need for a full unroll of the nonlinear iteration.

We start by reformulating our implicit scheme as an equation to find a root for:

\begin{equation}
    \label{eqn:rooteqn}
    g(y_{n}, y_{n+1}, \theta)=0.
\end{equation}

\noindent We solve this nonlinear equation for our future time point $y_{n+1}$ using a nonlinear solver such as Newton's method. Once we have the solution prediction \(y_{n+1}\) for the neural ODE at the time \(t_{n+1}\), we need to compute the gradient of this prediction with respect to the model parameters, \(\frac{\partial y_{n+1}}{\partial \theta}\), to update the neural network's parameters. While backpropagating through the multiple iterations of the nonlinear solver using automatic differentiation is possible, it is both resource-intensive and numerically unstable. Instead, we utilize the implicit function theorem for a more stable gradient calculation:

\begin{equation}
    \label{eqn:implicitfuntheorem}
    \frac{\partial y_{n+1}}{\partial \theta} = -\left(\frac{\partial g}{\partial y_{n+1}}\right)^{-1} \frac{\partial g}{\partial \theta}.
\end{equation}

\noindent The expression from the implicit function theorem is computed at the point where the nonlinear system is solved. If this solution lacks sufficient precision, the calculated gradients may suffer from numerical instability and inaccuracy.

In our previous work (see Ref.~\onlinecite{fronk2024trainingstiffneuralordinary}), we demonstrated successful training of stiff neural ODEs using the backward Euler, trapezoid, and Radau IIA methods. As a first-order method, backward Euler required an impractically high number of training data points and very small step sizes, underscoring the need for higher-order methods for training neural ODEs, especially when using adaptive step size solvers. Among the methods tested, the trapezoid method proved to be the most reliable, even though it is not the highest order method. Radau IIA, with its higher order, allowed for training with the largest step sizes, making it a strong candidate for handling stiff neural ODEs. Despite their effectiveness, these implicit methods are computationally expensive due to the need for solving nonlinear equations. Therefore, we aim to explore alternative approaches that avoid nonlinear solvers to reduce computational costs.

\subsection{Explicit Exponential Integration Methods}

Exponential integration methods \cite{trefethen2000spectral, boyd2001chebyshev, canuto1988spectral, ascher2008numerical, miranker2001numerical} approach the solution of an ODE by separating it into a linear component and the remaining nonlinear terms:

\begin{equation}
    \frac{dy}{dt} = f(t,y(t)) = L y + N(t,y)
    \label{eqn:linearizedODE}
\end{equation}

\noindent where the linear term, $L$, is the Jacobian evaluated at our initial condition $y_0$:

\begin{equation}
    L = \frac{df}{dy}(y_0)
\end{equation}

\noindent and the nonlinear term, $N$, represents the remaining components:

\begin{equation}
    N(t,y(t)) = f(t,y(t)) - Ly.
    \label{eqn:nonlinear}
\end{equation}

\noindent This decomposition allows for the exact integration of the linear component using the matrix exponential, which is particularly beneficial when the linear term represents the stiffest part of the dynamics, and the nonlinear part changes more slowly. By handling the stiff linear dynamics exactly and using an explicit method for the slower nonlinear part, we can better manage stiffness. However, this reformulation still leaves the question of how to effectively integrate the new ODE. 

\subsection*{Integrating Factor Methods}

Integrating factor methods \cite{lawson1967generalized, boyd2001chebyshev, canuto1988spectral, fornberg1999fast, trefethen2000spectral} use a change of variables

\begin{equation}
    w(t) = e^{-Lt} y(t),
\end{equation}

\noindent which, when differentiated with respect to $t$, results in:

\begin{equation}
    \frac{dw(t)}{dt} = e^{-Lt} \left (\frac{dy(t)}{dt} - Ly(t) \right ).
\end{equation}

\noindent By substituting Eqn. \ref{eqn:nonlinear} into our equation, we obtain:

\begin{equation}
    \frac{dw(t)}{dt} = e^{-Lt} N(t,y(t)) = e^{-Lt} N(t,e^{Lt}y(t)).
    \label{eqn:integratingfactor}
\end{equation}

\noindent Integrating factor methods are straightforward to derive because they allow the use of any standard numerical integration technique, such as multi-step or Runge-Kutta methods, to solve Eqn \ref{eqn:integratingfactor}. After applying the chosen integration method, the solution is converted back to the original variable $y$. Despite their simplicity, these methods are often discouraged in the literature because they can fail to compute fixed points correctly, resulting in greater errors compared to other exponential integration techniques. Despite the literature's discouragement, the integrating factor Euler method is the most robust exponential integration method we could find for neural ODEs. Applying the forward Euler method on Eqn \ref{eqn:integratingfactor} results in the A-stable first-order integrating factor Euler scheme \cite{lawson1967generalized, boyd2001chebyshev, canuto1988spectral, fornberg1999fast, trefethen2000spectral, cox2002exponential}, which has a local truncation error of $\frac{h^2}{2}LN$:

\begin{equation}
    y_{n+1} = e^{Lh} (  y_n + h N_n  ).
    \label{eqn:IFEuler}
\end{equation}

Despite an extensive search for a higher-order exponential integration method, our recent analysis (see Ref.~\onlinecite{fronk2024performanceevaluationsinglestepexplicit}) indicates that exponential integration methods fail to improve upon the first-order accuracy of IF Euler while remaining stable, revealing the IF Euler method as the only reliable choice for repeated, inexpensive integration in applications such as neural ODEs.  

\subsection*{Computing the Matrix Exponential}

\noindent Computing the matrix exponential \cite{10.1137/1020098, doi:10.1137/S00361445024180, del2003survey, RUIZ2016370, ARIOLI1996111}, \( e^A \), for these exponential integration schemes can be computationally expensive. For dense matrices, a common approach is the scaling and squaring method with a Pade approximation, which involves scaling down the matrix, applying a rational approximation, and squaring the result. This method has a complexity of \( \mathcal{O}(n^3) \), making it impractical for large matrices. If \( A \) is diagonalizable, \( e^A \) can be computed as \( e^A = V e^\Lambda V^{-1} \), where \( \Lambda \) is diagonal, but this still has a complexity of \( \mathcal{O}(n^3) \) due to costly matrix multiplications and inversions. For large sparse matrices, Krylov subspace methods, such as Arnoldi iteration, can reduce the complexity to \( \mathcal{O}(n^2) \). However, these iterative methods may require numerous iterations to reach the desired accuracy, which can increase computational costs and reduce stability. As the dimension of the ODE system grows, the cost of computing the matrix exponential, which scales as \( \mathcal{O}(n^2) \) in the best case scenario for iterative methods, can become prohibitively expensive.

\FloatBarrier
\clearpage
\section{Results}

We start by testing our methodology on the stiff univariate linear equation (Example 1), which serves as a standard benchmark for stability analysis in different ODE discretization methods, making it an ideal initial baseline. Next, we conduct tests on a 10D stiff linear system (Example 2) to verify whether the performance characteristics observed earlier for the univariate linear equation persist as the system dimension is increased. A 3D nonlinear stiff system of ODEs, hand-designed for this test, is used in Example 3 to examine the methods' capabilities in nonlinear scenarios. We conclude with Example 4, which tests the methodology on a stiff Van der Pol oscillator, a model where classical implicit methods fail during neural ODE training. This example showcases the stability of the exponential integrating factor Euler method, which succeeds where classical implicit methods fall short. We produced training data for each of these model and analyzed how effectively the different single-step integration methods could recover the ODE model.

\subsection{Example 1: Stiff Linear Model}

In our first test, we consider the stiff linear equation:

\begin{equation}
    \frac{dy}{dt} = -10000y, \quad y(0)=1000, \quad t \in [0,0.01].
\end{equation}

\begin{figure*}
    \centering
    \includegraphics[width=0.9\linewidth]{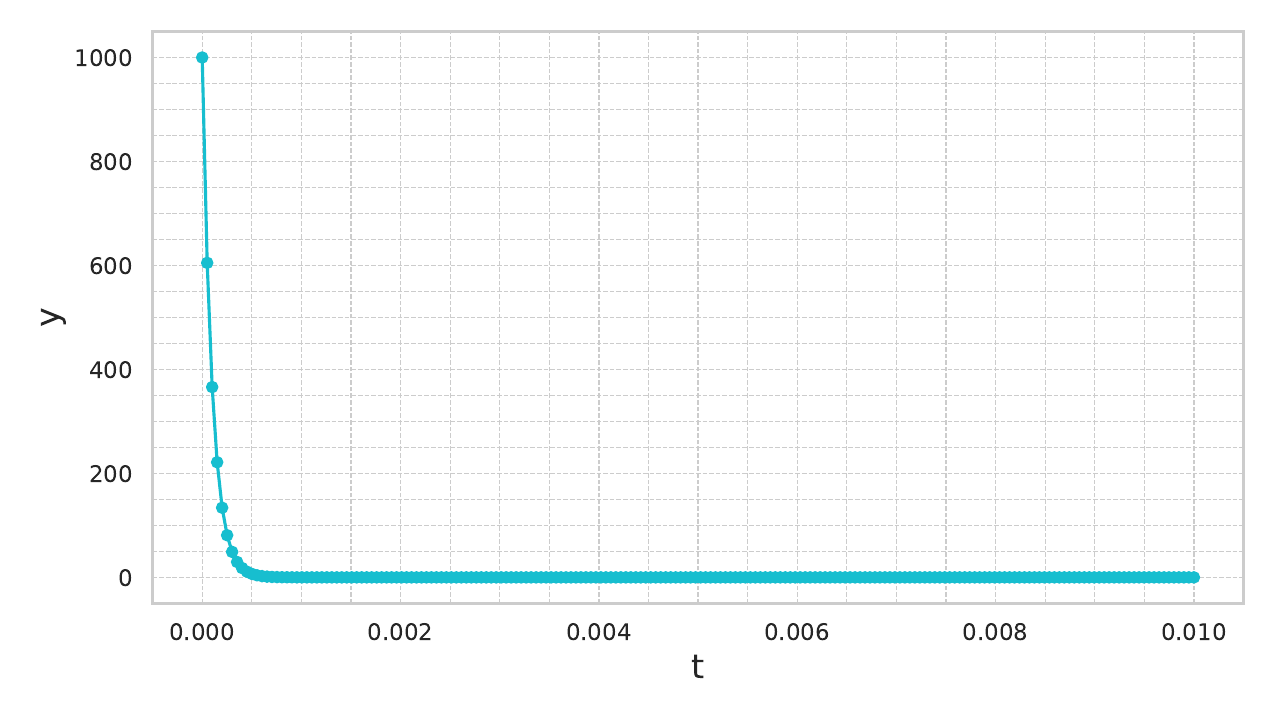}
    \caption{For the equation $y'=-10000y$, we plot the training dataset consisting of 200 data points uniformly distributed across the time interval.}
    \label{fig:1D_Training_Data}
\end{figure*}

\noindent The model depicted in Figure \ref{fig:1D_Training_Data} is a standard benchmark for evaluating the stability of numerical methods on stiff systems of ODEs. Its high degree of stiffness, driven by the large negative coefficient, results in an initial sharp transient followed by a slower decay to equilibrium. The slow decay region is the stiff region where explicit methods struggle due to stability limitations.

We used the SciPy \cite{2020SciPy-NMeth} Radau solver to generate training data. Figure \ref{fig:1D_Training_Data} displays the training data, consisting of the case with 100 uniformly spaced time points. As discussed in the methods section, we used the discretize-then-optimize strategy, segmenting the data into $n-1$ time intervals, each treated as an IVP between adjacent time points. During each epoch, our custom JAX \cite{jax2018github, deepmind2020jax} implementation of the solvers discussed in the methods section solves these $n-1$ IVPs in parallel.

\begin{figure*}
    \centering
    \includegraphics[width=0.9\linewidth]{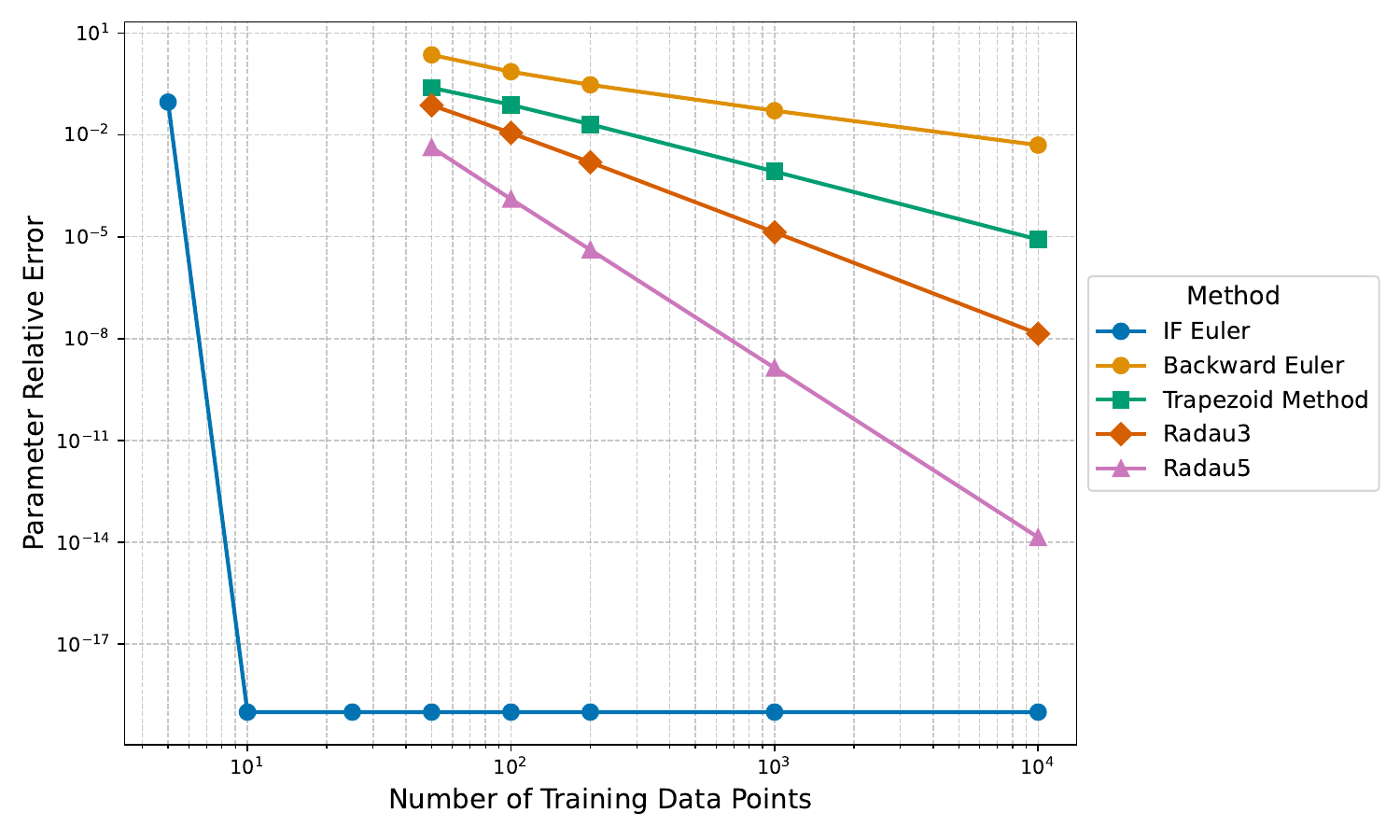}
    \caption{The fractional parameter relative error (non-percentage) is plotted against the number of training data points for the equation $y'=-10000y$. For comparison, we show the explicit exponential integrating factor Euler method method alongside a few implicit schemes.}
    \label{fig:Linear_Error_vs_n}
\end{figure*}

To establish a preliminary benchmark, we tested the implicit methods backward Euler, trapezoid, Radau3, and Radau5, using training data consisting of 50 to 10,000 uniformly spaced points in time. We then tested the performance of the exponential integrating factor Euler method on the same training data. For each experiment, we recovered the equation learned by the $\pi$-net V1 polynomial neural network and evaluated its accuracy by calculating the fractional relative error for each of the model's inferred parameters. Table 1 shows the recovered equations for varying amounts of training data, while Figure \ref{fig:Linear_Error_vs_n} plots the fractional relative error versus the number of data points. Remarkably, the exponential integrating factor Euler method accurately recovered the exact equation for all datasets, down to numerical precision. To determine its limitations, we gradually reduced $n$ until the method failed, which occurred with just 5 training points. This finding suggests that the method's precise integration of linear ODEs makes it phenomenally well-suited for training linear ODEs.

\begin{table*}\centering
\ra{1.3}
\begin{tabular}{@{}rrrrcrrrcrrr@{}}\toprule
& \multicolumn{3}{c}{$y'=-10000y$} \\
\cmidrule{2-4} 
& $n$ & & Equation Learned \\ \midrule
Exponential IF Euler Method\\
& $5$ & $y'$ =& $  - 10927.463178714 y  +  3.42071479417 \cdot 10^{-5}   $ \\
& $10$ & $y'$ =& $  - 9999.999999999999999 y - 1.49710130728912 \cdot 10^{-13}   $ \\
& $25$ & $y'$ =& $  - 9999.999999999999999 y + 1.02210460890402 \cdot 10^{-11}   $ \\
& $50$ & $y'$ =& $  - 9999.999999999999999 y + 8.74976237336876 \cdot 10^{-12}   $ \\
& $100$ & $y'$ =& $  - 9999.999999999999999 y - 2.00875522420467 \cdot 10^{-12}  $ \\
& $200$ & $y'$ =& $  - 9999.999999999999999 y + 1.70268949846918 \cdot 10^{-11}   $ \\
& $1000$ & $y'$ =& $  - 9999.999999999999999 y - 4.88937945290105 \cdot 10^{-12}   $ \\
& $10000$ & $y'$ =&  $  - 9999.999999999999999 y - 1.54427227469678 \cdot 10^{-11}   $ \\
Backward Euler\\
& $50$ & $y'$ =& $- 32814.7600708 y + 8.44726018571 \cdot 10^{-11} $ \\
& $100$ & $y'$ =& $- 17284.1957884 y - 6.65078880626 \cdot 10^{-12}$ \\
& $200$ & $y'$ =& $- 12992.0930002 y + 9.10992647337 \cdot 10^{-12} $ \\
& $1000$ & $y'$ =& $- 10517.6269781 y + 1.65232316322 \cdot 10^{-11}$ \\
& $10000$ & $y'$ =& $- 10050.1721181 y + 2.23615888785 \cdot 10^{-11}$ \\
Trapezoid Method\\
& $50$ & $y'$ =& $- 7546.32076209 y - 4.25736335374 \cdot 10^{-12}$ \\
& $100$ & $y'$ =& $- 9228.38069787 y - 2.92238168648 \cdot 10^{-11}$ \\
& $200$ & $y'$ =& $- 9794.7490243 y + 1.16207823362 \cdot 10^{-8} $ \\
& $1000$ & $y'$ =& $- 9991.65833324 y + 4.60703162317 \cdot 10^{-12} $ \\
& $10000$ & $y'$ =& $- 9999.91665083 y - 4.69854345928 \cdot 10^{-11}$ \\
Radau3\\
& $50$ & $y'$ =& $- 9251.48316114 y - 1.0589518914 \cdot 10^{-11}$ \\
& $100$ & $y'$ =& $- 9885.79527641 y - 6.93292052555 \cdot 10^{-12}$ \\
& $200$ & $y'$ =& $- 9984.36246911 y - 5.90851257248 \cdot 10^{-11}$ \\
& $1000$ & $y'$ =& $- 9999.86426589 y + 8.25230499684 \cdot 10^{-12}$ \\
& $10000$ & $y'$ =& $- 9999.99986144 y + 4.41857407072 \cdot 10^{-12}$ \\
Radau5\\
& $50$ & $y'$ =& $- 10042.9715925 y - 5.81410253423 \cdot 10^{-11}$ \\
& $100$ & $y'$ =& $- 10001.2886455 y + 8.15234805477 \cdot 10^{-12}$ \\
& $200$ & $y'$ =& $- 10000.0413085 y - 2.43678163057 \cdot 10^{-11}$ \\
& $1000$ & $y'$ =& $- 10000.0000137 y - 1.01896661931 \cdot 10^{-11}$ \\
& $10000$ & $y'$ =& $- 10000.00000000014 y + 1.06509360321 \cdot 10^{-11} $ \\
\bottomrule
\end{tabular}
\caption{Comparison of recovered equations for $y'=-10000y$ using various stiff ODE methods. 
 The results of the explicit exponential integrating factor Euler method are compared with selected implicit single-step methods, with a varying number of training points ($n$) equally spaced in time over the time interval.}
\end{table*}

\clearpage
\FloatBarrier
\subsection{Example 2: 10-Dimensional Stiff Linear Model}

Our next example is a 10-dimensional stiff linear model:

\begin{equation}
\label{eqn:10D-Linear}
\begin{aligned}
    &\frac{dy_0}{dt} = -10 y_{0} + 5 y_{1}, \\
    &\frac{dy_1}{dt} = 5 y_{0} - 20 y_{1} + 5 y_{2}, \\
    &\frac{dy_2}{dt} = 5 y_{1} - 50 y_{2} + 5 y_{3}, \\
    &\frac{dy_3}{dt} = 5 y_{2} - 100 y_{3} + 5 y_{4}, \\
    &\frac{dy_4}{dt} = 5 y_{3} - 500 y_{4} + 5 y_{5}, \\
    &\frac{dy_5}{dt} = 5 y_{4} - 1000 y_{5} + 5 y_{6}, \\
    &\frac{dy_6}{dt} = 5 y_{5} - 5000 y_{6} + 5 y_{7}, \\
    &\frac{dy_7}{dt} = 5 y_{6} - 10000 y_{7} + 5 y_{8}, \\
    &\frac{dy_8}{dt} = 5 y_{7} - 20000 y_{8} + 5 y_{9}, \\
    &\frac{dy_9}{dt} = 5 y_{8} - 50000 y_{9}, \\
    &y_0(0) = 20, \quad y_1(0) = 20, \quad y_2(0) = 20, \quad y_3(0) = 20, \\
     &y_4(0) = 20, \quad y_5(0) = 20, \quad y_6(0) = 20, \quad y_7(0) = 20, \\
      &y_8(0) = 20, \quad y_9(0) = 20,  \quad t \in [0, 0.4]. 
\end{aligned}
\end{equation}

\begin{figure*}
    \centering
    \includegraphics[width=0.9\linewidth]{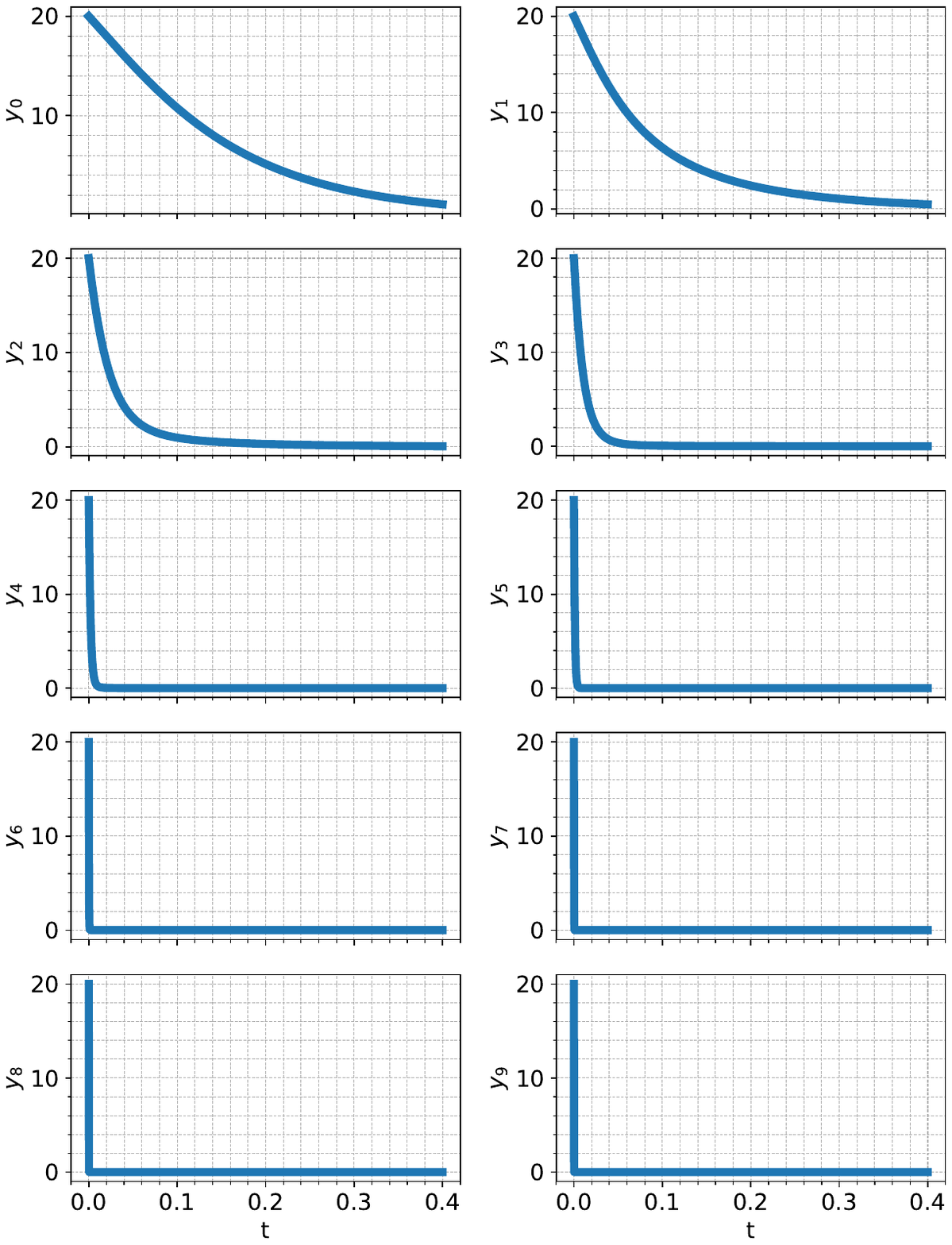}
    \caption{The training data for the 10-dimensional linear stiff system of ODEs (Example 2).  The stiffness can be seen by the fact that  each of the variables approach steady-state at a different timescale.}
    \label{fig:10D_training_data}
\end{figure*}

\begin{figure*}
    \centering
    \includegraphics[angle=90, width=0.9\linewidth]{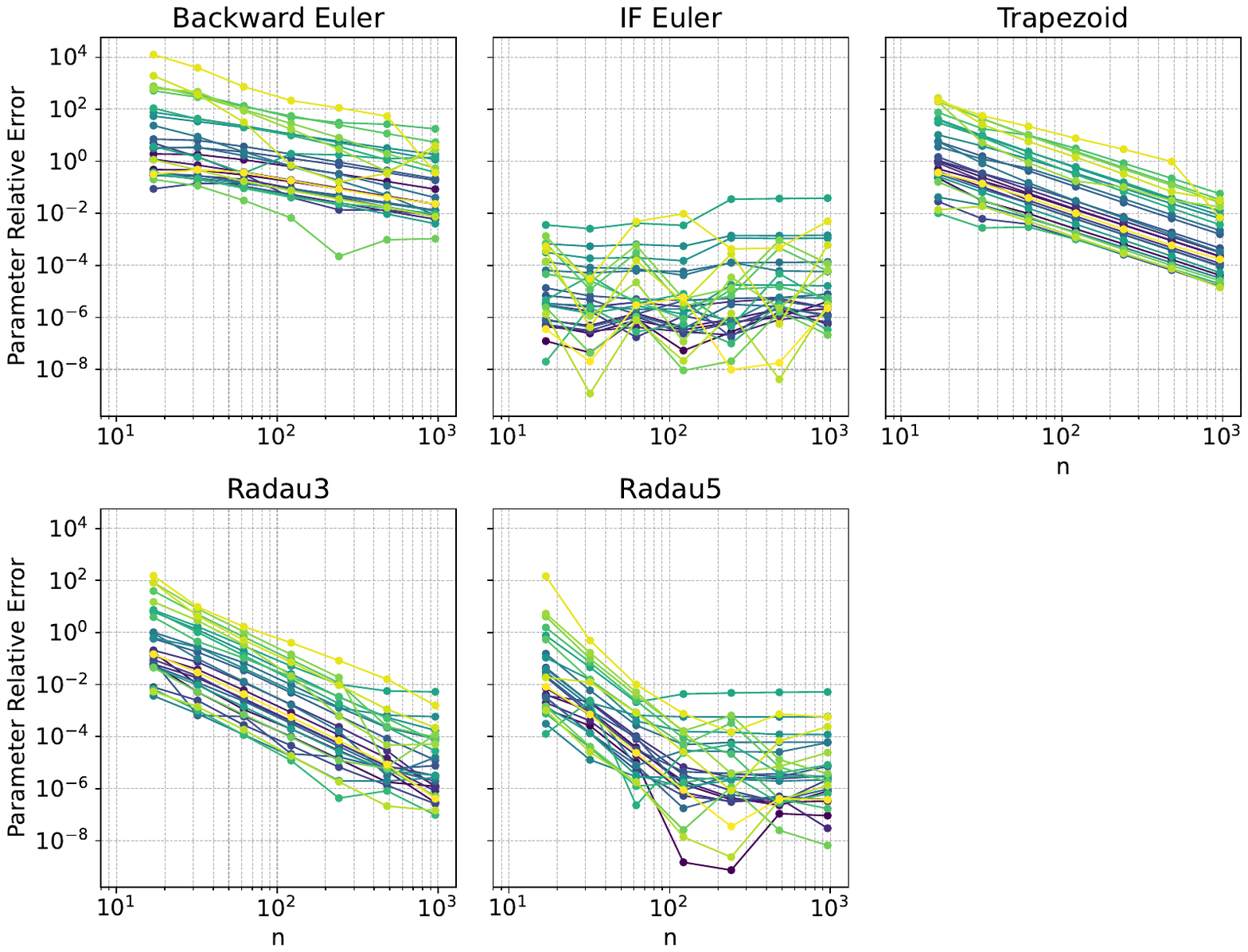}
    \caption{For Example 2, we plot the fractional parameter relative error (not a percentage) against the number of training data points for the exponential integrating factor Euler, Backward Euler, Trapezoid, Radau3, and Radau5 integration schemes. Subfigures are organized by integration scheme, showing the 28 parameters of the linear model, with each parameter represented by a unique color.}
    \label{fig:10D_Linear_Error_vs_n}
\end{figure*}

\noindent The intriguing results of the 1-dimensional stiff linear model led us to test a higher-dimensional version to see if the observed patterns would be present with more variables. We engineered this toy model to have a concise parameter count of 28, considering its larger 10-dimensional size, while also ensuring that it had a high degree of stiffness. As shown in Figure \ref{fig:10D_training_data}, the model’s 10 variables exhibit varying time scales, contributing to its stiffness. We generated training data using the same integration technique and discretize-then-optimize approach outlined in Example 1. We tested the same integration schemes as before for varying numbers of training data points ranging from 10 to 1000 time points. For each trial, we trained the polynomial neural network, recovered the equation learned by the polynomial neural network, and then assessed the accuracy of the identified parameters by calculating the fractional relative error. Tables 2 through 6 display the equations recovered with a limited dataset of $n=17$ using the various stiff integration schemes. Figure \ref{fig:10D_Linear_Error_vs_n} shows the fractional relative error vs the number of training data points. The findings from this experiment are consistent with those seen in the 1D stiff linear system. The IF Euler method demonstrates remarkable accuracy in reproducing the correct equations, attributable to its precise time integration for linear systems.  The IF Euler method excels in accuracy when training data is scarce. As the training data increases, it matches the performance of the Radau5 method.

% Table of coefficients for example 2 - Backward Euler
\begin{turnpage}
\begin{table*}\centering
\ra{1.3}
\begin{tabular}{@{}r r r r r r r r r r r r@{}}\toprule
& $y_0$ & $y_1$ & $y_2$ & $y_3$ & $y_4$ & $y_5$ & $y_6$ & $y_7$ & $y_8$ & $y_9$ & $b$ \\
\cmidrule{2-12}
$y_0'$ & -14.7836 & 14.5781 & -8.4924 & 6.1276 & -7.1492 & 9.2093 & -22.9949 & 43.7244 & -33.0434 & -4.6647 & -13.4089 \\
$y_1'$ & -0.9625 & -13.2098 & 3.5689 & 2.1554 & -8.7490 & 13.5899 & -24.3976 & 21.5091 & 9.5306 & -25.6484 & -4.2411 \\
$y_2'$ & 12.2460 & -20.1279 & -45.6400 & 21.9133 & -39.8748 & 49.9039 & -25.5398 & 0.9724 & -9.5407 & -25.4553 & -3.9481 \\
$y_3'$ & 4.4018 & -2.5930 & -30.7243 & -64.8650 & -10.5595 & 26.6175 & -49.4778 & 50.3420 & 4.6271 & -18.8469 & -2.9835 \\
$y_4'$ & 31.3069 & -64.2247 & 100.1548 & -113.9335 & -648.0725 & 277.7614 & -73.1184 & -34.0929 & -18.1808 & -132.0993 & -15.0393 \\
$y_5'$ & -1.5240 & 2.1398 & -2.9124 & 13.1950 & -381.7649 & -619.1673 & -13.6554 & -2.1618 & -17.9343 & -15.1688 & 3.1873 \\
$y_6'$ & -11.7862 & 13.4005 & 15.8586 & -32.6563 & 232.0548 & -544.3964 & -6743.0863 & 3716.9921 & -3362.2256 & -57.0894 & -35.5675 \\
$y_7'$ & 2.4226 & -4.9786 & -11.4434 & 46.9336 & -87.0531 & 96.8002 & -2565.2063 & -7986.4847 & -3870.5058 & 265.9424 & -43.0494 \\
$y_8'$ & -240.2358 & 136.8660 & 200.3869 & 93.8648 & -0.3288 & -579.9877 & -364.3709 & 3100.9704 & -42918.6704 & 9706.6464 & -29.5044 \\
$y_9'$ & 47.5369 & 32.2822 & -12.0926 & -144.5704 & -107.4558 & -103.5888 & 2507.9738 & 8568.2294 & -62821.6534 & -34034.1024 & 6.1518 \\
\bottomrule
\end{tabular}
\caption{Coefficient matrix for Example 2 using the backward Euler method, with training data of size $n=17$. Rows correspond to equations for $y_0'$ to $y_9'$, and columns represent coefficients for $y_0$ to $y_9$ and $b$.}
\end{table*}
\end{turnpage}

% Table of coefficients for example 2 (IF Euler)
\begin{turnpage}
\begin{table*}\centering
\ra{1.3}
\begin{tabular}{@{}r r r r r r r r r r r r@{}}\toprule
& $y_0$ & $y_1$ & $y_2$ & $y_3$ & $y_4$ & $y_5$ & $y_6$ & $y_7$ & $y_8$ & $y_9$ & $b$ \\
\cmidrule{2-12}
$y_0'$ & -10.00000012 & 5.00000255 & -4.9835e-06 & 8.1128e-06 & -5.9082e-05 & 1.6135e-04 & -1.9551e-03 & 7.2901e-03 & -1.2028e-02 & 9.9142e-03 & 6.3076e-07 \\
$y_1'$ & 4.99999616 & -19.99999127 & 4.99998262 & 2.7227e-05 & -1.8137e-04 & 4.8569e-04 & -5.7973e-03 & 2.1693e-02 & -3.6013e-02 & 2.9876e-02 & 1.3840e-06 \\
$y_2'$ & -6.8013e-06 & 5.00001494 & -50.00002647 & 5.00003462 & -1.7539e-04 & 4.4168e-04 & -4.9090e-03 & 1.8096e-02 & -2.9763e-02 & 2.4560e-02 & 2.7067e-06 \\
$y_3'$ & 1.9772e-05 & -4.2538e-05 & 5.00006759 & -100.000083 & 5.00031826 & -7.5264e-04 & 7.9394e-03 & -2.9512e-02 & 4.9415e-02 & -4.1535e-02 & -9.0289e-06 \\
$y_4'$ & 2.1494e-04 & -4.4515e-04 & 6.4912e-04 & 4.99930709 & -499.9983232 & 4.99663035 & 2.6865e-02 & -9.4728e-02 & 1.5490e-01 & -1.2943e-01 & -1.1058e-04 \\
$y_5'$ & 2.7239e-04 & -5.4942e-04 & 7.5903e-04 & -7.6034e-04 & 5.00156986 & -1000.0028415 & 5.01787658 & -5.9551e-02 & 9.3807e-02 & -7.7032e-02 & -1.4969e-04 \\
$y_6'$ & 6.8083e-06 & -1.0720e-05 & 6.0775e-06 & -1.5962e-07 & -1.4257e-05 & 5.00002322 & -5000.00009875 & 5.00023081 & -3.1307e-04 & 2.3873e-04 & -6.7437e-06 \\
$y_7'$ & 1.8040e-05 & 1.2622e-05 & 2.2877e-05 & 1.8154e-05 & 5.4586e-05 & 1.1048e-04 & 5.00344242 & -9999.97583443 & 5.00666282 & 1.2957e-03 & 4.4818e-06 \\
$y_8'$ & 2.3246e-05 & -6.3319e-06 & 4.3669e-06 & 3.0982e-05 & -6.0586e-05 & 1.2727e-04 & 1.4412e-03 & 5.00070203 & -19999.97289518 & 5.00179466 & 8.0888e-06 \\
$y_9'$ & -7.1509e-04 & 1.2919e-03 & -1.2730e-03 & 8.6851e-04 & -4.7200e-04 & 3.6830e-04 & 3.7458e-05 & 3.1387e-04 & 5.00270677 & -49999.98230818 & 5.1730e-04 \\
\bottomrule
\end{tabular}
\caption{Coefficient matrix for Example 2 using the exponential integrating factor Euler method, with training data of size $n=17$. Rows correspond to equations for $y_0'$ to $y_9'$, and columns represent coefficients for $y_0$ to $y_9$ and $b$.}
\end{table*}
\end{turnpage}

% Table of coefficients for example 2 (Trapezoid method)
\begin{turnpage}
\begin{table*}\centering
\ra{1.3}
\begin{tabular}{@{}r r r r r r r r r r r r@{}}\toprule
& $y_0$ & $y_1$ & $y_2$ & $y_3$ & $y_4$ & $y_5$ & $y_6$ & $y_7$ & $y_8$ & $y_9$ & $b$ \\
\cmidrule{2-12}
$y_0'$ & -12.2841 & 7.5659 & -1.3450 & 0.8123 & -0.7411 & 0.8725 & -1.3970 & 1.6440 & -0.0788 & -0.6361 & 8.9738 \\
$y_1'$ & 10.0457 & -28.0456 & 10.3687 & -2.8062 & 1.6789 & -1.8580 & 2.4636 & -2.1082 & -0.7982 & 1.2405 & -1.1987 \\
$y_2'$ & -2.1192 & 9.2792 & -48.5786 & -2.3307 & 11.1623 & -8.2038 & -2.7017 & 1.6107 & 2.9435 & 1.2452 & 0.7548 \\
$y_3'$ & 1.8465 & -7.9471 & 33.8552 & -132.5776 & 23.1214 & -2.0256 & -7.9026 & -4.1798 & -0.4616 & 0.4157 & 1.7303 \\
$y_4'$ & -10.5070 & 20.1116 & -27.1287 & 33.8439 & -479.0886 & -46.8417 & 13.1621 & 11.7461 & 4.3970 & 1.3456 & 9.0838 \\
$y_5'$ & 10.7318 & -25.3258 & 45.5948 & -55.2749 & 214.7758 & -1260.9753 & 158.0829 & -15.8488 & -83.3956 & -3.2163 & -1.7247 \\
$y_6'$ & 8.1982 & -10.2751 & -7.4064 & 24.6267 & -107.0063 & 203.2613 & -4949.6236 & -373.8089 & 251.1799 & 32.3822 & -6.5969 \\
$y_7'$ & -9.1156 & 14.3421 & -0.1697 & -16.9750 & 86.2549 & -150.8973 & 1051.5011 & -11612.7722 & 995.1486 & -42.8935 & -1.2609 \\
$y_8'$ & -3.4403 & 5.7688 & -1.6838 & -5.2228 & 30.6558 & -45.2697 & 9.5634 & 992.9547 & -20273.6674 & 1371.9336 & 2.0202 \\
$y_9'$ & 6.4000 & -7.1302 & -1.7645 & 6.6440 & -6.2281 & -7.2155 & 179.3156 & -623.2880 & 1073.6613 & -31710.0728 & -15.5599 \\
\bottomrule
\end{tabular}
\caption{Coefficient matrix for Example 2 using the trapezoid method, with training data of size $n=17$. Rows correspond to equations for $y_0'$ to $y_9'$, and columns represent coefficients for $y_0$ to $y_9$ and $b$.}
\end{table*}
\end{turnpage}

\begin{turnpage}
\begin{table*}\centering
\ra{1.3}
\begin{tabular}{@{}r r r r r r r r r r r r@{}}\toprule
& $y_0$ & $y_1$ & $y_2$ & $y_3$ & $y_4$ & $y_5$ & $y_6$ & $y_7$ & $y_8$ & $y_9$ & $b$ \\
\cmidrule{2-12}
$y_0'$ & -10.6184 & 5.8565 & -0.6182 & 0.4069 & -0.3637 & 0.4038 & -0.3397 & 0.0666 & 0.1203 & 0.1399 & 1.4865 \\
$y_1'$ & 5.6450 & -21.2199 & 6.0343 & -0.6120 & 0.3982 & -0.4303 & 0.5109 & -0.3395 & -0.2678 & 0.2866 & 0.6341 \\
$y_2'$ & -0.0308 & 5.4551 & -50.3977 & 4.6113 & 1.6938 & -2.0441 & 1.0361 & 0.5374 & -0.6088 & -0.5580 & -0.1280 \\
$y_3'$ & -0.8256 & 1.0006 & 7.8841 & -105.1974 & 10.1152 & -2.0665 & -2.7752 & -0.2631 & 0.7860 & 0.8779 & 0.4902 \\
$y_4'$ & -0.6272 & 1.4739 & -2.9491 & 9.7142 & -501.8329 & 1.9133 & 4.8428 & -0.8331 & -0.8517 & -1.6286 & 0.2075 \\
$y_5'$ & 0.1299 & -0.3729 & 0.9980 & -2.0665 & 36.0745 & -1046.1492 & 41.5037 & -3.6098 & -24.3847 & -1.3311 & -0.0092 \\
$y_6'$ & 0.0535 & 0.2841 & -1.7797 & 3.1163 & -14.0448 & 39.2398 & -5033.1376 & -14.4872 & 62.0771 & -34.8870 & -0.2965 \\
$y_7'$ & -0.1360 & 0.1221 & 0.3451 & -0.8056 & 3.6929 & -7.9307 & 202.7018 & -10444.8240 & 422.1632 & -162.2079 & 0.1911 \\
$y_8'$ & -0.1298 & 0.3306 & -0.2036 & -0.5346 & 6.4674 & -13.4232 & 70.5477 & 80.2254 & -20106.8303 & 416.0239 & -0.2106 \\
$y_9'$ & 0.5998 & -0.2810 & -0.4414 & -1.6470 & 5.3721 & 6.8277 & -64.2909 & -98.7092 & 763.4865 & -42656.8993 & -1.2393 \\
\bottomrule
\end{tabular}
\caption{Coefficient matrix for Example 2 using the Radau3 method, with training data of size $n=17$. Rows correspond to equations for $y_0'$ to $y_9'$, and columns represent coefficients for $y_0$ to $y_9$ and $b$.}
\end{table*}
\end{turnpage}

% Table of coefficients for example 2 (Radau5)
\begin{turnpage}
\begin{table*}\centering
\ra{1.3}
\begin{tabular}{@{}r r r r r r r r r r r r@{}}\toprule
& $y_0$ & $y_1$ & $y_2$ & $y_3$ & $y_4$ & $y_5$ & $y_6$ & $y_7$ & $y_8$ & $y_9$ & $b$ \\
\cmidrule{2-12}
$y_0'$ & -9.9886 & 4.9797 & 1.9113e-02 & -1.3712e-02 & 1.4039e-02 & -1.7570e-02 & 3.9613e-02 & -7.1463e-02 & 5.3269e-02 & -1.4328e-02 & -1.1274e-03 \\
$y_1'$ & 5.0368 & -20.0400 & 5.0185 & -1.0287e-02 & 9.4517e-03 & -1.1820e-02 & 2.4294e-02 & -3.7282e-02 & 1.5168e-02 & 9.8373e-03 & -1.4774e-01 \\
$y_2'$ & -1.3782e-01 & 5.2228 & -50.1614 & 5.0971 & -7.9934e-02 & 9.7715e-02 & -2.1924e-01 & 4.0289e-01 & -3.1798e-01 & 1.0683e-01 & 2.1727e-02 \\
$y_3'$ & 4.2241e-02 & -8.3983e-02 & 4.8936 & -99.6840 & 4.2357 & 1.0302 & -2.3242 & 3.8436 & -2.0194 & -4.8388e-01 & -8.3587e-03 \\
$y_4'$ & -5.3687e-02 & 1.2765e-01 & -2.4722e-01 & 5.1493 & -499.8454 & 4.8000 & 3.4440e-01 & -6.4920e-01 & 6.0098e-01 & -3.0590e-01 & 1.6773e-02 \\
$y_5'$ & 6.3468e-02 & -1.4571e-01 & 2.6941e-01 & -3.4961e-01 & 4.4616 & -998.5717 & 1.1826 & 6.9396 & -5.6358 & 2.0738 & -2.1570e-02 \\
$y_6'$ & -2.9514e-02 & 6.6207e-02 & -1.1522e-01 & 1.4387e-01 & -6.2377e-01 & 5.0006 & -4996.1840 & -2.7608 & 7.8833 & -4.9776 & 1.0654e-02 \\
$y_7'$ & 6.9226e-02 & -1.5295e-01 & 2.5519e-01 & -2.9377e-01 & 8.2757e-01 & -1.6728 & 2.3061 & -9985.7412 & -15.5597 & 15.1246 & -2.5701e-02 \\
$y_8'$ & -3.3101e-02 & 8.7182e-02 & -1.9103e-01 & 2.5379e-01 & -7.9914e-01 & 1.4719 & -13.9334 & 31.5335 & -20020.2009 & 5.0941 & 2.8611e-03 \\
$y_9'$ & 5.0149e-01 & -9.1170e-01 & 4.8384e-01 & 7.3468e-01 & -9.1965 & 19.0498 & -122.4756 & 408.5092 & -724.1976 & -50411.6780 & -1.5342e-01 \\
\bottomrule
\end{tabular}
\caption{Coefficient matrix for Example 2 using the Radau5 method, with training data of size $n=17$. Rows correspond to equations for $y_0'$ to $y_9'$, and columns represent coefficients for $y_0$ to $y_9$ and $b$.}
\end{table*}
\end{turnpage}

\clearpage
\FloatBarrier
\subsection{Example 3: 3D Nonlinear Model}

In our third example, we consider the following stiff nonlinear system of ODEs:

\begin{equation}
\label{eqn:example3}
\begin{aligned}
    &\frac{dy_1}{dt} = -500 y_1 + 3.8 y_2^2 + 1.35 y_3,   \\
    &\frac{dy_2}{dt} = 0.82 y_1 - 24 y_2 + 7.5 y_3^2,  \\
    &\frac{dy_3}{dt} = -0.5 y_1^2 + 1.85 y_2 - 6.5 y_3^2,  \\
    &y_1(0) = 15, \quad y_2(0) = 7, \quad y_3(0) = 10, \quad t \in [0, 5].
\end{aligned}
\end{equation}

\begin{figure*}
    \centering
    \includegraphics[width=0.9\linewidth]{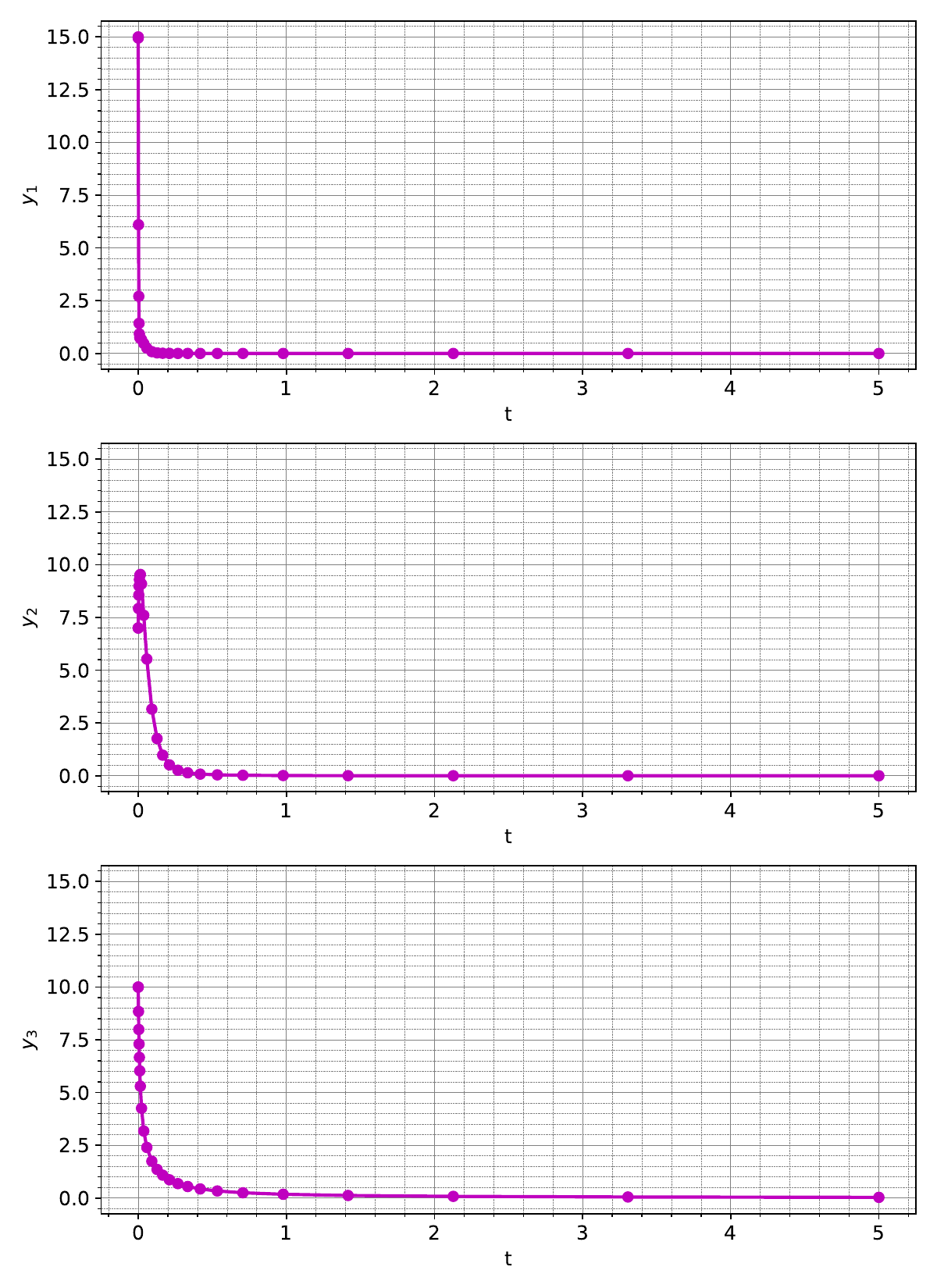}
    \caption{Training data for Example 3, a 3D nonlinear stiff ODE system.  This is the training data for the model corresponding to $n=25$ data points.}
    \label{fig:3D_Nonlinear_Training_data_subfigures}
\end{figure*}

\noindent Figure \ref{fig:3D_Nonlinear_Training_data_subfigures} shows the training region for this model. We specifically created this stiff neural ODE toy problem to be three-dimensional and contain only linear and quadratic terms. This increases the difficulty of numerical integration slightly but remains a small enough model for a comprehensive analysis of the model's nine identified parameters.

\begin{figure*}
    \centering
    \includegraphics[width=0.90\linewidth]{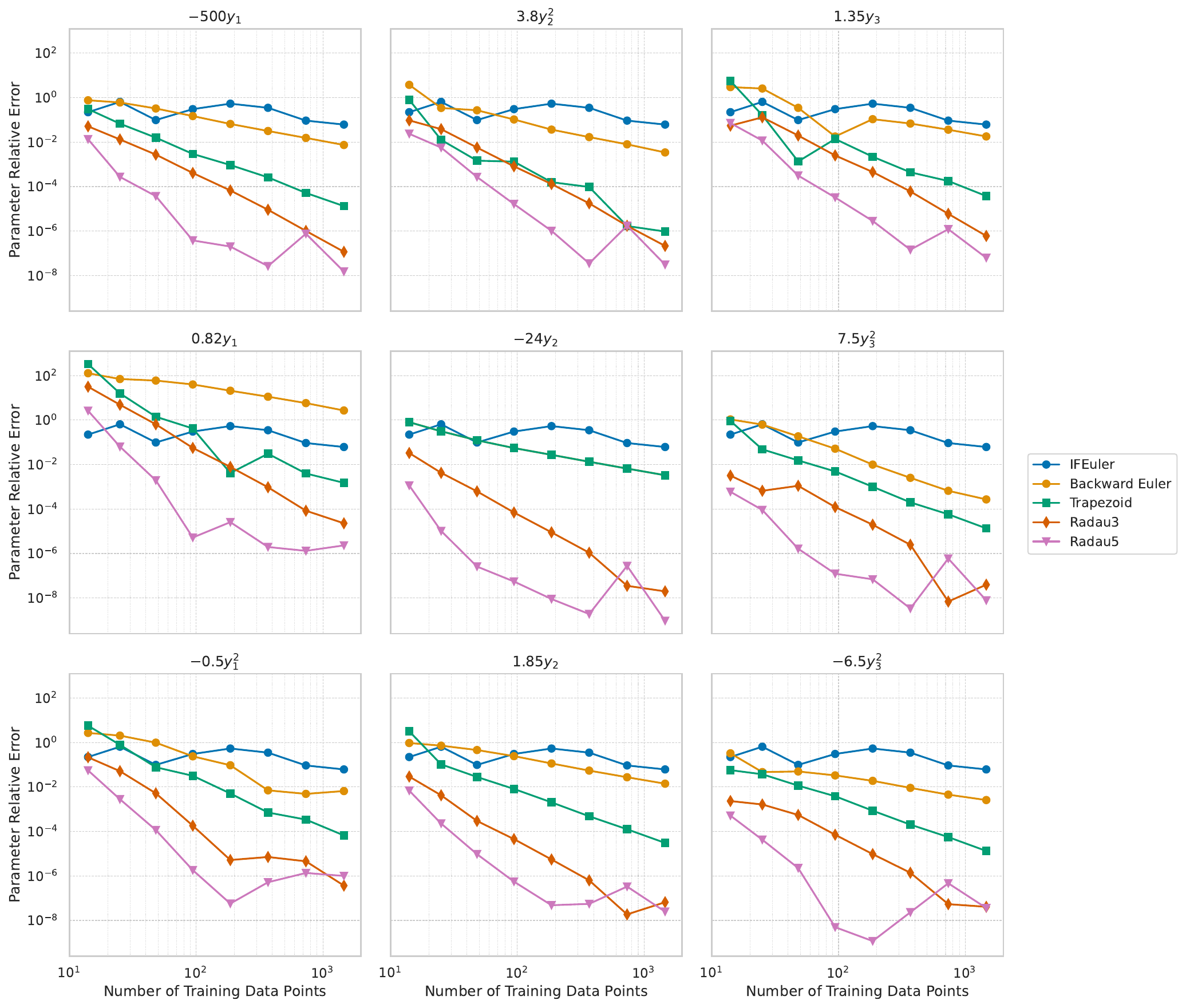}
    \caption{For Example 3, we plot the fractional parameter relative error (not a percentage) against the number of training data points for the exponential integrating factor Euler, Backward Euler, Trapezoid, Radau3, and Radau5 integration schemes. The symbols displayed in the legend represent each integration method. Each subfigure corresponds to a specific term in the ODE, as indicated by the title of each subfigure.}
    \label{fig:Error_vs_n_3D_nonlinear}
\end{figure*}

We generated training data using the same integration technique and discretize-then-optimize approach outlined in Example 1. We varied $n$, the total number of data points, to study how different amounts of training data affect the accuracy of the inferred parameters. We trained the polynomial neural ODE on each dataset, recovered the corresponding equations from the polynomial neural ODE, and assessed the accuracy of the inferred parameters by computing their fractional relative error. Tables 7, 8, 9, and 10 show the recovered equations for datasets consisting of 1467, 369, 94, and 48 data points respectively, while Figure \ref{fig:Error_vs_n_3D_nonlinear} illustrates how the fractional relative error changes with the amount of training data. Although the exponential integrating factor Euler method is more cost-effective, the tables and figures highlight that it achieves the lowest overall performance for this model.

%table 1 of example 2 - n = 1467
\begin{turnpage}
\begin{table*}\centering
\ra{1.3}
\begin{tabular}{@{}rrrrcrrrcrrr@{}}\toprule
& \multicolumn{3}{c}{ 
    $\begin{aligned}
    y_1' &= -500 y_1 + 3.8 y_2^2 + 1.35 y_3, \\
    y_2' &= 0.82 y_1 - 24 y_2 + 7.5 y_3^2, \\
    y_3' &= -0.5 y_1^2 + 1.85 y_2 - 6.5 y_3^2
    \end{aligned}$
    } \\
\cmidrule{2-4} 
& $n$ & & Equation Learned \\ \midrule

Exponential IF Euler\\
& $1467$ & $y_1'$ =& $1.51968131012 y_{1}^{2} + 7.46491448656 y_{1} y_{2} - 4.4360187548 y_{1} y_{3} - 530.76798674 y_{1} + 3.23377291781 y_{2}^{2}   $\\ 
& & & $+ 0.743119055549 y_{2} y_{3} - 0.748806101868 y_{2} - 0.100968562079 y_{3}^{2} + 1.5824687907 y_{3} - 0.0132703872778$  \\
&    & $y_2'$ =&  $ - 0.0147620261445 y_{1}^{2} - 0.112013689897 y_{1} y_{2} + 0.195326314213 y_{1} y_{3} - 0.552493633928 y_{1} + 0.0388233505295 y_{2}^{2} $\\ 
&&& $- 0.11936884918 y_{2} y_{3} - 24.0551857426 y_{2} + 7.615389906 y_{3}^{2} - 0.113933526112 y_{3} - 0.00208780383352$ \\
&    & $y_3'$ =&  $ - 0.565017339582 y_{1}^{2} - 0.0226566807316 y_{1} y_{2} + 0.386791942974 y_{1} y_{3} - 2.60233668838 y_{1} + 0.017698859856 y_{2}^{2} $\\ 
&&& $- 0.0113901872367 y_{2} y_{3} + 1.91229571307 y_{2} - 6.50755324083 y_{3}^{2} + 0.0212679661193 y_{3} - 8.05782938895 \cdot 10^{-5}$ \\

Backward Euler\\
& $1467$ & $y_1'$ =& $0.0384064527302 y_{1}^{2} + 0.204017196661 y_{1} y_{2} - 0.184501300028 y_{1} y_{3} - 503.719351023 y_{1} + 3.81306747361 y_{2}^{2} $\\ 
& & & $+ 0.030695277063 y_{2} y_{3} + 0.0138909801005 y_{2} - 0.0389371667798 y_{3}^{2} + 1.37431473228 y_{3} - 0.000981504457688$  \\
&    & $y_2'$ =& $0.0200787332271 y_{1}^{2} + 0.284950328409 y_{1} y_{2} + 0.00680359696552 y_{1} y_{3} - 1.3914960677 y_{1} + 0.0120642998079 y_{2}^{2} $ \\ 
&&& $- 0.00393399669949 y_{2} y_{3} - 24.0794531237 y_{2} + 7.49797005252 y_{3}^{2} + 0.0118933425261 y_{3} - 0.000244780267862 $ \\
&    & $y_3'$ =& $- 0.503277229186 y_{1}^{2} - 0.0936986490342 y_{1} y_{2} - 0.0504958796953 y_{1} y_{3} + 0.969202527143 y_{1} - 0.000853834054426 y_{2}^{2} $ \\ 
&&& $- 0.00916411500017 y_{2} y_{3} + 1.82402554075 y_{2} - 6.48322554439 y_{3}^{2} - 0.0124328526534 y_{3} + 0.000363525556554 $ \\

Trapezoid Method\\
& $1467$ & $y_1'$ =& $- 0.000214354771517 y_{1}^{2} - 0.000974387910813 y_{1} y_{2} + 0.00118871099339 y_{1} y_{3} - 499.993355185 y_{1} + 3.79999636531 y_{2}^{2} $\\ 
& & & $- 7.6389755885 \cdot 10^{-5} y_{2} y_{3} - 6.49292906827 \cdot 10^{-5} y_{2} + 0.00010856272938 y_{3}^{2} + 1.34994917202 y_{3} + 1.7778522614 \cdot 10^{-6}$  \\
&    & $y_2'$ =& $- 9.40804457738 \cdot 10^{-5} y_{1}^{2} - 0.000738401443738 y_{1} y_{2} + 0.000464830560082 y_{1} y_{3} + 0.821234592939 y_{1} - 4.98824711904 \cdot 10^{-6} y_{2}^{2} $ \\ 
&&& $- 4.45829770482 \cdot 10^{-5} y_{2} y_{3} - 23.9998826649 y_{2} + 7.50010148624 y_{3}^{2} - 3.27328943583 \cdot 10^{-5} y_{3} + 8.82781486984 \cdot 10^{-7} $ \\
&    & $y_3'$ =& $- 0.499966869555 y_{1}^{2} + 0.000213011256539 y_{1} y_{2} + 2.40008585951 \cdot 10^{-5} y_{1} y_{3} - 0.00124987034708 y_{1} - 2.79306256779 \cdot 10^{-6} y_{2}^{2}  $\\ 
&&& $+ 4.49476766209 \cdot 10^{-5} y_{2} y_{3} + 1.85005644455 y_{2} - 6.50008628126 y_{3}^{2} + 4.10531838766 \cdot 10^{-5} y_{3} - 1.29639482701 \cdot 10^{-6} $ \\

Radau3\\
& $1467$ & $y_1'$ =& $- 1.58565658931 \cdot 10^{-6} y_{1}^{2} - 9.0220800355 \cdot 10^{-6} y_{1} y_{2} + 5.34341576966 \cdot 10^{-6} y_{1} y_{3} - 499.999941385 y_{1} + 3.79999918096 y_{2}^{2} $ \\ 
& & & $+ 2.82082509581 \cdot 10^{-6} y_{2} y_{3} - 6.09956239295 \cdot 10^{-7} y_{2} - 2.19977308688 \cdot 10^{-6} y_{3}^{2} + 1.35000081713 y_{3} - 6.6422572823 \cdot 10^{-8}$  \\
&    & $y_2'$ =& $- 4.10343797756 \cdot 10^{-7} y_{1}^{2} - 7.61410752424 \cdot 10^{-7} y_{1} y_{2} + 2.85277919953 \cdot 10^{-6} y_{1} y_{3} + 0.819981280026 y_{1} + 1.08825064073 \cdot 10^{-7} y_{2}^{2} $ \\ 
&&& $- 3.29411976356 \cdot 10^{-7} y_{2} y_{3} - 23.9999995232 y_{2} + 7.50000029772 y_{3}^{2} + 1.89809742021 \cdot 10^{-8} y_{3} + 2.28812493605 \cdot 10^{-9} $ \\
&    & $y_3'$ =& $- 0.499999816321 y_{1}^{2} + 1.05053334608 \cdot 10^{-6} y_{1} y_{2} - 1.03391317172 \cdot 10^{-6} y_{1} y_{3} + 1.56427731904 \cdot 10^{-6} y_{1} - 9.58120432784 \cdot 10^{-8} y_{2}^{2}  $\\ 
&&& $+ 2.22479763452 \cdot 10^{-7} y_{2} y_{3} + 1.85000012191 y_{2} - 6.50000026651 y_{3}^{2} + 6.09264676245 \cdot 10^{-8} y_{3} - 2.12359008156 \cdot 10^{-9} $ \\

\bottomrule
\end{tabular}
\caption{Comparison of recovered equations for Example 3 using various stiff ODE methods. The results of the exponential integrating factor Euler method are compared with selected implicit single-step methods for $n=1467$ training data points.}
\end{table*}
\end{turnpage}

%table 2 of Example 2
\begin{turnpage}
\begin{table*}\centering
\ra{1.3}
\begin{tabular}{@{}rrrrcrrrcrrr@{}}\toprule
& \multicolumn{3}{c}{ 
    $\begin{aligned}
    y_1' &= -500 y_1 + 3.8 y_2^2 + 1.35 y_3, \\
    y_2' &= 0.82 y_1 - 24 y_2 + 7.5 y_3^2, \\
    y_3' &= -0.5 y_1^2 + 1.85 y_2 - 6.5 y_3^2
    \end{aligned}$
    } \\
\cmidrule{2-4} 
& $n$ & & Equation Learned \\ \midrule
Exponential IF Euler\\
& $369$ & $y_1'$ =& $ 6.14645119889 y_{1}^{2} + 16.2082700784 y_{1} y_{2} - 35.9173670546 y_{1} y_{3} - 325.227702362 y_{1} - 0.251338737648 y_{2}^{2}   $\\ 
& & & $ + 9.27315536001 y_{2} y_{3} - 6.86945893334 y_{2} - 6.01411434056 y_{3}^{2} + 2.67692884315 y_{3} + 5.31091112729 $  \\
&    & $y_2'$ =&  $ 0.166538587879 y_{1}^{2} + 0.298078666792 y_{1} y_{2} - 0.64485493486 y_{1} y_{3} + 1.73022817427 y_{1} + 0.0177848633454 y_{2}^{2} $\\ 
&&& $  - 0.193436525603 y_{2} y_{3} - 24.1352320224 y_{2} + 7.76315657869 y_{3}^{2} - 0.382945373609 y_{3} - 0.028553067734  $ \\
&    & $y_3'$ =&  $ - 0.597760821756 y_{1}^{2} - 0.104976548032 y_{1} y_{2} + 0.587460505627 y_{1} y_{3} - 3.00013788369 y_{1} + 0.00854108649362 y_{2}^{2} $\\ 
&&& $+ 0.08255522856 y_{2} y_{3} + 1.94693472977 y_{2} - 6.61314573829 y_{3}^{2} + 0.081432123729 y_{3} + 0.0473882095799$ \\

Backward Euler\\
& $369$ & $y_1'$ =& $0.127904202437 y_{1}^{2} + 0.797524296352 y_{1} y_{2} - 0.618853828345 y_{1} y_{3} - 515.835938291 y_{1} + 3.86400911308 y_{2}^{2} $\\ 
& & & $ + 0.100596632061 y_{2} y_{3} + 0.0702514281895 y_{2} - 0.144792011507 y_{3}^{2} + 1.44253523851 y_{3} - 0.00338565729602$  \\
&    & $y_2'$ =& $0.0543209089864 y_{1}^{2} + 1.04579692729 y_{1} y_{2} + 0.160430638587 y_{1} y_{3} - 8.37766824113 y_{1} + 0.0654796412269 y_{2}^{2}  $\\ 
&&& $- 0.052630703619 y_{2} y_{3} - 24.3209270038 y_{2} + 7.51898512114 y_{3}^{2} + 0.0406798119831 y_{3} - 0.000613776339809 $ \\
&    & $y_3'$ =& $- 0.503529354636 y_{1}^{2} - 0.329057521838 y_{1} y_{2} - 0.245153506381 y_{1} y_{3} + 3.85641129319 y_{1} - 0.00849774225915 y_{2}^{2}  $\\ 
&&& $- 0.026326319065 y_{2} y_{3} + 1.74932608891 y_{2} - 6.44088493316 y_{3}^{2} - 0.048392901916 y_{3} + 0.00140276670564 $ \\

Trapezoid Method\\
& $369$ & $y_1'$ =& $- 0.00268209780416 y_{1}^{2} - 0.0144329545699 y_{1} y_{2} + 0.0150266836908 y_{1} y_{3} - 499.869671739 y_{1} + 3.79963430099 y_{2}^{2} $\\ 
& & & $- 0.000501143645822 y_{2} y_{3} - 0.00134685668212 y_{2} + 0.00131751622372 y_{3}^{2} + 1.3494057239 y_{3} + 1.90847684559 \cdot 10^{-5}$  \\
&    & $y_2'$ =& $- 0.00114528900208 y_{1}^{2} - 0.0106000717948 y_{1} y_{2} + 0.00558187371416 y_{1} y_{3} + 0.84505905353 y_{1} - 0.000230229973003 y_{2}^{2} $ \\ 
&&& $- 0.000472976562712 y_{2} y_{3} - 23.9980217168 y_{2} + 7.50149654194 y_{3}^{2} - 0.000524988587427 y_{3} + 4.37738386438 \cdot 10^{-5} $ \\
&    & $y_3'$ =& $- 0.499644737787 y_{1}^{2} + 0.00285307226783 y_{1} y_{2} + 0.00129595454858 y_{1} y_{3} - 0.0229524685876 y_{1} + 2.98567414682 \cdot 10^{-5} y_{2}^{2} $ \\ 
&&& $+ 0.000590463049014 y_{2} y_{3} + 1.85089037745 y_{2} - 6.50132012361 y_{3}^{2} + 0.00064671327907 y_{3} - 2.36121005357 \cdot 10^{-5} $ \\

Radau3\\
& $369$ & $y_1'$ =& $- 7.61974773614 \cdot 10^{-5} y_{1}^{2} - 0.000498327709512 y_{1} y_{2} + 0.000176699490206 y_{1} y_{3} - 499.995450213 y_{1} + 3.79993286416 y_{2}^{2} $\\ 
& & & $+ 0.000217509049623 y_{2} y_{3} - 7.01297138472 \cdot 10^{-5} y_{2} - 0.00017272660185 y_{3}^{2} + 1.35008196827 y_{3} - 7.79201917524 \cdot 10^{-6}$  \\
&    & $y_2'$ =& $- 1.12137698434 \cdot 10^{-5} y_{1}^{2} - 2.36452610708 \cdot 10^{-5} y_{1} y_{2} + 0.000100980310608 y_{1} y_{3} + 0.819221372929 y_{1} + 2.4231234723 \cdot 10^{-6} y_{2}^{2}  $\\ 
&&& $- 1.52392601319 \cdot 10^{-5} y_{2} y_{3} - 23.9999739349 y_{2} + 7.50001854938 y_{3}^{2} + 1.28276717158 \cdot 10^{-7} y_{3} + 1.54151826501 \cdot 10^{-7} $ \\
&    & $y_3'$ =& $- 0.499996445415 y_{1}^{2} + 8.11323903976 \cdot 10^{-6} y_{1} y_{2} - 2.60051361773 \cdot 10^{-5} y_{1} y_{3} + 0.000206368555709 y_{1} - 1.47477468822 \cdot 10^{-6} y_{2}^{2}  $\\ 
&&& $+ 6.83199869368 \cdot 10^{-6} y_{2} y_{3} + 1.84999882977 y_{2} - 6.50000902966 y_{3}^{2} + 3.08390275437 \cdot 10^{-6} y_{3} - 1.29968263207 \cdot 10^{-7} $ \\

\bottomrule
\end{tabular}
\caption{Comparison of recovered equations for Example 3 using various stiff ODE methods. The results of the exponential integrating factor Euler method are compared with selected implicit single-step methods for $n=369$ training data points.}
\end{table*}
\end{turnpage}

%table 3 of Example 2
\begin{turnpage}
\begin{table*}\centering
\ra{1.3}
\begin{tabular}{@{}rrrrcrrrcrrr@{}}\toprule
& \multicolumn{3}{c}{ 
    $\begin{aligned}
    y_1' &= -500 y_1 + 3.8 y_2^2 + 1.35 y_3, \\
    y_2' &= 0.82 y_1 - 24 y_2 + 7.5 y_3^2, \\
    y_3' &= -0.5 y_1^2 + 1.85 y_2 - 6.5 y_3^2
    \end{aligned}$
    } \\
\cmidrule{2-4} 
& $n$ & & Equation Learned \\ \midrule
Exponential IF Euler\\
& $94$ & $y_1'$ =& $ - 1.78703554807 y_{1}^{2} + 1.15294499924 y_{1} y_{2} + 15.1750052006 y_{1} y_{3} - 652.313344767 y_{1} + 1.1639973465 y_{2}^{2}   $\\ 
& & & $ + 3.30628446599 y_{2} y_{3} - 3.42081114919 y_{2} - 3.23945008745 y_{3}^{2} + 4.96697425016 y_{3} + 211.76388665 $  \\
&    & $y_2'$ =&  $ 0.435541301972 y_{1}^{2} - 1.82821183452 y_{1} y_{2} - 1.83416277298 y_{1} y_{3} + 20.4914706073 y_{1} + 0.30499722096 y_{2}^{2} $\\ 
&&& $  - 1.27277464229 y_{2} y_{3} - 24.3846298027 y_{2} + 9.00824418215 y_{3}^{2} - 1.11852463442 y_{3} - 7.17082844967  $ \\
&    & $y_3'$ =&  $ - 0.270455465228 y_{1}^{2} + 2.06150354604 y_{1} y_{2} - 1.05996016603 y_{1} y_{3} - 4.94243827667 y_{1} - 0.074093652478 y_{2}^{2} $\\ 
&&& $  + 0.334558514439 y_{2} y_{3} + 1.73194413238 y_{2} - 6.8929220562 y_{3}^{2} + 0.717153752746 y_{3} + 1.63417705184   $ \\

Backward Euler\\
& $94$ & $y_1'$ =& $0.628429562731 y_{1}^{2} + 3.83347856227 y_{1} y_{2} - 2.44967509189 y_{1} y_{3} - 574.356843937 y_{1} + 4.19704419353 y_{2}^{2} $\\ 
& & & $ - 0.0507037997737 y_{2} y_{3} + 0.639626711303 y_{2} - 0.173536682816 y_{3}^{2} + 1.37419372243 y_{3} + 0.029408755562$  \\
&    & $y_2'$ =& $0.309621055396 y_{1}^{2} + 4.03977576762 y_{1} y_{2} + 0.137940498732 y_{1} y_{3} - 31.8965120786 y_{1} + 0.34747681308 y_{2}^{2}  $\\ 
&&& $- 0.523025515187 y_{2} y_{3} - 25.3243973178 y_{2} + 7.88766384537 y_{3}^{2} + 0.00379946914356 y_{3} + 0.00975557096823 $ \\
&    & $y_3'$ =&  $- 0.619343143732 y_{1}^{2} - 1.68702386966 y_{1} y_{2} - 0.365857735248 y_{1} y_{3} + 13.2374339888 y_{1} - 0.0081667697632 y_{2}^{2}  $\\ 
&&& $- 0.090220830018 y_{2} y_{3} + 1.39932809859 y_{2} - 6.2844789663 y_{3}^{2} - 0.168276206228 y_{3} + 0.0038905402588 $ \\

Trapezoid Method\\
& $94$ & $y_1'$ =& $- 0.0588646265565 y_{1}^{2} - 0.27012838604 y_{1} y_{2} + 0.345013331708 y_{1} y_{3} - 498.52567995 y_{1} + 3.80504335312 y_{2}^{2} $\\ 
& & & $- 0.0336776819517 y_{2} y_{3} - 0.0127434045842 y_{2} + 0.0397814247098 y_{3}^{2} + 1.33125963682 y_{3} + 0.000982572835436$  \\
&    & $y_2'$ =& $- 0.0450668815737 y_{1}^{2} - 0.23297656984 y_{1} y_{2} + 0.244365221953 y_{1} y_{3} + 0.476040309876 y_{1} + 0.00869338753968 y_{2}^{2}  $\\ 
&&& $- 0.0299121574977 y_{2} y_{3} - 23.9612707912 y_{2} + 7.53681178165 y_{3}^{2} - 0.0129193788573 y_{3} + 0.000718185923707 $ \\
&    & $y_3'$ =& $- 0.48410748181 y_{1}^{2} + 0.0805474244863 y_{1} y_{2} - 0.0356373166309 y_{1} y_{3} - 0.188145335786 y_{1} - 0.00387491084267 y_{2}^{2}  $\\ 
&&& $+ 0.0159521980287 y_{2} y_{3} + 1.86508272955 y_{2} - 6.52489715182 y_{3}^{2} + 0.0109209453638 y_{3} - 0.000391914744422 $ \\

Radau3\\
& $94$ & $y_1'$ =& $ - 0.00365909221556 y_{1}^{2} - 0.0223368391164 y_{1} y_{2} + 0.0104314867655 y_{1} y_{3} - 499.795116702 y_{1} + 3.79689360939 y_{2}^{2} $\\ 
& & & $ + 0.00967601990744 y_{2} y_{3} - 0.00291843193932 y_{2} - 0.00759128716958 y_{3}^{2} + 1.3534044362 y_{3} - 0.000302176679636$  \\
&    & $y_2'$ =& $ - 0.000535842537473 y_{1}^{2} - 0.000757268154964 y_{1} y_{2} + 0.00536390419245 y_{1} y_{3} + 0.774043217656 y_{1} + 3.64111573268 \cdot 10^{-5} y_{2}^{2} $ \\ 
&&& $ - 0.000681620496574 y_{2} y_{3} - 23.9983199818 y_{2} + 7.50092201196 y_{3}^{2} + 9.65683744741 \cdot 10^{-6} y_{3} + 6.23666981312 \cdot 10^{-6} $ \\
&    & $y_3'$ =& $ - 0.500090837014 y_{1}^{2} - 0.000436889453858 y_{1} y_{2} - 3.30120008436 \cdot 10^{-5} y_{1} y_{3} + 0.0078554942054 y_{1} + 2.05865827988 \cdot 10^{-5} y_{2}^{2}  $\\ 
&&& $ + 0.000252612907303 y_{2} y_{3} + 1.8499162119 y_{2} - 6.50046305391 y_{3}^{2} + 0.000191341009964 y_{3} - 7.60979876269 \cdot 10^{-6} $ \\

\bottomrule
\end{tabular}
\caption{Comparison of recovered equations for Example 3 using various stiff ODE methods. The results of the exponential integrating factor Euler method are compared with selected implicit single-step methods for $n=94$ training data points.}
\end{table*}
\end{turnpage}

%table 4 of Example 3
\begin{turnpage}
\begin{table*}\centering
\ra{1.3}
\begin{tabular}{@{}rrrrcrrrcrrr@{}}\toprule
& \multicolumn{3}{c}{ 
    $\begin{aligned}
    y_1' &= -500 y_1 + 3.8 y_2^2 + 1.35 y_3, \\
    y_2' &= 0.82 y_1 - 24 y_2 + 7.5 y_3^2, \\
    y_3' &= -0.5 y_1^2 + 1.85 y_2 - 6.5 y_3^2
    \end{aligned}$
    } \\
\cmidrule{2-4} 
& $n$ & & Equation Learned \\ \midrule
Exponential IF Euler\\
& $48$ & $y_1'$ =& $ 2.25022013897 y_{1}^{2} + 12.6446422113 y_{1} y_{2} - 14.09885659 y_{1} y_{3} - 450.853573172 y_{1} - 1.1212500601 y_{2}^{2}   $\\ 
& & & $ + 8.93619930997 y_{2} y_{3} - 5.49047776771 y_{2} - 8.50372012866 y_{3}^{2} + 10.2653249402 y_{3} + 208.15129309 $  \\
&    & $y_2'$ =&  $ - 0.150872114012 y_{1}^{2} - 11.088465431 y_{1} y_{2} + 0.96684003193 y_{1} y_{3} + 60.9862908595 y_{1} + 1.2571939785 y_{2}^{2} $\\ 
&&& $   - 4.15898664239 y_{2} y_{3} - 21.1548016982 y_{2} + 11.6977197696 y_{3}^{2} - 3.36676106746 y_{3} - 31.0843168357  $ \\
&    & $y_3'$ =&  $ - 0.255798609357 y_{1}^{2} + 3.31471715978 y_{1} y_{2} - 1.1294891592 y_{1} y_{3} - 12.0853834847 y_{1} - 0.157468068538 y_{2}^{2} $\\ 
&&& $  + 0.714791019829 y_{2} y_{3} + 1.53275092969 y_{2} - 7.27455175132 y_{3}^{2} + 1.22512284862 y_{3} + 5.78628817503   $ \\

Backward Euler\\
& $48$ & $y_1'$ =& $ 1.82910880564 y_{1}^{2} + 8.77522528086 y_{1} y_{2} - 6.26048259683 y_{1} y_{3} - 662.79250133 y_{1} + 4.82816206628 y_{2}^{2} $\\ 
& & & $ - 0.815104217543 y_{2} y_{3} + 1.67112795975 y_{2} + 0.37076347215 y_{3}^{2} + 0.880167900117 y_{3} + 0.12015687347$  \\
&    & $y_2'$ =& $ 0.821363180957 y_{1}^{2} + 6.93094406635 y_{1} y_{2} - 0.948084348591 y_{1} y_{3} - 48.1915891864 y_{1} + 0.863938300142 y_{2}^{2}  $\\ 
&&& $ - 1.60765034538 y_{2} y_{3} - 26.9635740125 y_{2} + 8.89368929416 y_{3}^{2} - 0.289151173692 y_{3} + 0.0411175673398 $ \\
&    & $y_3'$ =& $ - 0.991176710831 y_{1}^{2} - 3.56810157187 y_{1} y_{2} + 0.754691504145 y_{1} y_{3} + 16.2729284742 y_{1} + 0.030668469293 y_{2}^{2} $ \\ 
&&& $ - 0.142692804358 y_{2} y_{3} + 1.0023097408 y_{2} - 6.17589039304 y_{3}^{2} - 0.280743480112 y_{3} + 0.00520012969419 $ \\

Trapezoid Method\\
& $48$ & $y_1'$ =& $ - 0.241197044867 y_{1}^{2} - 1.20003110545 y_{1} y_{2} + 1.27269890609 y_{1} y_{3} - 492.058348031 y_{1} + 3.79442941946 y_{2}^{2} $\\ 
& & & $ - 0.0365651879257 y_{2} y_{3} - 0.134839497026 y_{2} + 0.0728196465699 y_{3}^{2} + 1.35182703869 y_{3} - 0.00547754968617 $  \\
&    & $y_2'$  =& $ - 0.173412349544 y_{1}^{2} - 0.896381372281 y_{1} y_{2} + 0.93067994678 y_{1} y_{3} - 0.339973344039 y_{1} + 0.0227551702672 y_{2}^{2}  $\\ 
&&& $ - 0.0841998788607 y_{2} y_{3} - 23.8484059563 y_{2} + 7.61494110899 y_{3}^{2} - 0.0396939701894 y_{3} + 0.00205575932009 $ \\
&    & $y_3'$ =& $ - 0.46146772343 y_{1}^{2} + 0.221411627156 y_{1} y_{2} - 0.0322504040161 y_{1} y_{3} - 0.884561939734 y_{1} - 0.00166472898244 y_{2}^{2} $ \\ 
&&& $ + 0.0341049503204 y_{2} y_{3} + 1.90210507435 y_{2} - 6.57477886428 y_{3}^{2} + 0.0386700013798 y_{3} - 0.00136341520966 $ \\

Radau3\\
& $48$ & $y_1'$ =& $ - 0.00365909221556 y_{1}^{2} - 0.0223368391164 y_{1} y_{2} + 0.0104314867655 y_{1} y_{3} - 499.795116702 y_{1} + 3.79689360939 y_{2}^{2}  $\\ 
& & & $ + 0.00967601990744 y_{2} y_{3} - 0.00291843193932 y_{2} - 0.00759128716958 y_{3}^{2} + 1.3534044362 y_{3} - 0.000302176679636 $  \\
&    & $y_2'$ =& $ - 0.000535842537473 y_{1}^{2} - 0.000757268154964 y_{1} y_{2} + 0.00536390419245 y_{1} y_{3} + 0.774043217656 y_{1} + 3.64111573268 \cdot 10^{-5} y_{2}^{2}  $\\ 
&&& $ - 0.000681620496574 y_{2} y_{3} - 23.9983199818 y_{2} + 7.50092201196 y_{3}^{2} + 9.65683744741 \cdot 10^{-6} y_{3} + 6.23666981312 \cdot 10^{-6} $ \\
&    & $y_3'$ =& $ - 0.500090837014 y_{1}^{2} - 0.000436889453858 y_{1} y_{2} - 3.30120008436 \cdot 10^{-5} y_{1} y_{3} + 0.0078554942054 y_{1} + 2.05865827988 \cdot 10^{-5} y_{2}^{2}  $\\ 
&&& $ + 0.000252612907303 y_{2} y_{3} + 1.8499162119 y_{2} - 6.50046305391 y_{3}^{2} + 0.000191341009964 y_{3} - 7.60979876269 \cdot 10^{-6} $ \\

\bottomrule
\end{tabular}
\caption{Comparison of recovered equations for Example 3 using various stiff ODE methods. The results of the exponential integrating factor Euler method are compared with selected implicit single-step methods for $n=48$ training data points.}
\end{table*}
\end{turnpage}

\clearpage
\subsection{Example 4}

\begin{equation}
\label{eqn:example3}
\begin{aligned}
    &\frac{dx}{dt} = y,   \\
    &\frac{dy}{dt} = 1000 y - 1000 x^2 y - x,  \\
    &x(0) = 1, \quad y(0)=0, \quad t \in [0,1300]
\end{aligned}
\end{equation}

\begin{figure*}
    \centering
    \includegraphics[width=0.9\linewidth]{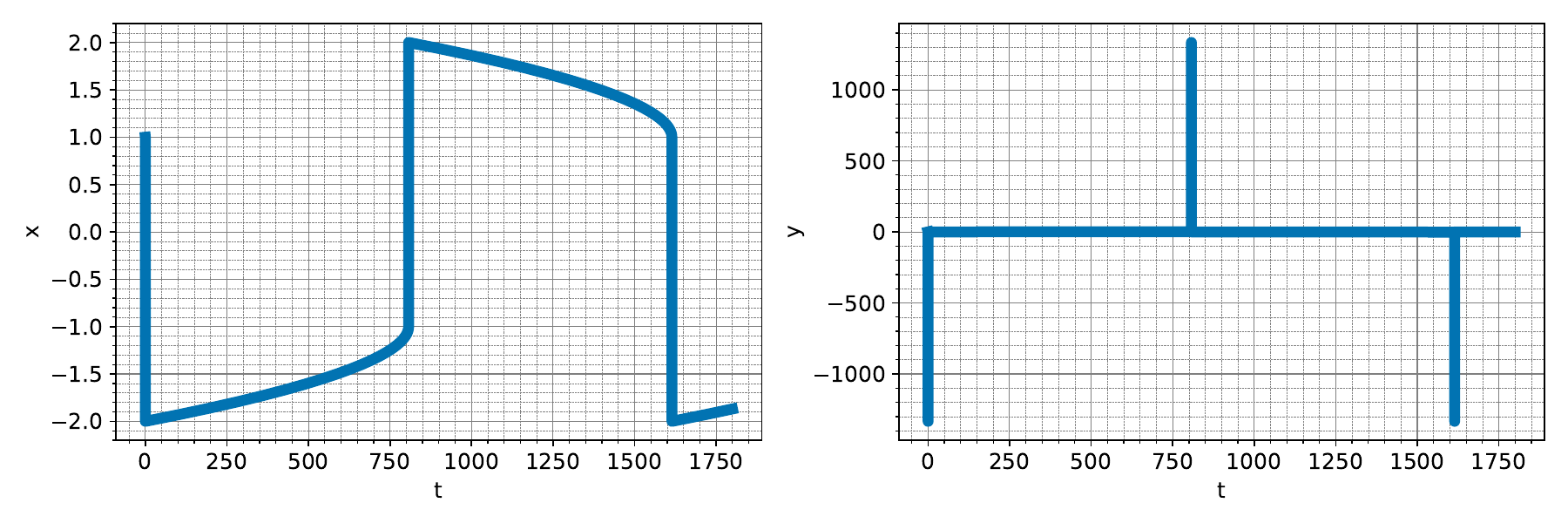}
    \caption{The training region for the stiff Van der Pol model (see Eqn.\ref{eqn:example3})}
    \label{fig:Training_data_vanderpol}
\end{figure*}

\noindent As a nonconservative system with nonlinear damping, the Van der Pol model is a widely used tool for studying relaxation oscillations \cite{VanderPolModel}. Its applications span across multiple fields, including neuroscience, seismology, and speech analysis, modeling phenomena such as neuron activity, fault line tremors, and vocal cord vibrations \cite{FITZHUGH1961445, 4066548, doi:10.1142/S0218127499001620, doi:10.1121/1.4798467}. The model’s stiffness can be adjusted through the parameter $\mu$, with higher values making the system stiffer, ideal for testing neural ODEs designed for stiff problems. In our study, we set $\mu=1000$, creating an extremely stiff ODE system. The method for generating training data remained unchanged from the process outlined in Example 1. Figure \ref{fig:Training_data_vanderpol} shows the training data, revealing rapid transients with large spikes in values followed by stiff regions where changes are minimal. Every implicit integration technique we tried failed to properly train the stiff neural ODE on this complex system, with convergence issues disrupting the process entirely. Even when we provided training data with more frequent time intervals to improve stability, the increase in computational demand made this approach unfeasible and did not solve the underlying convergence problem.

\begin{figure*}
    \centering
    \includegraphics[width=0.9\linewidth]{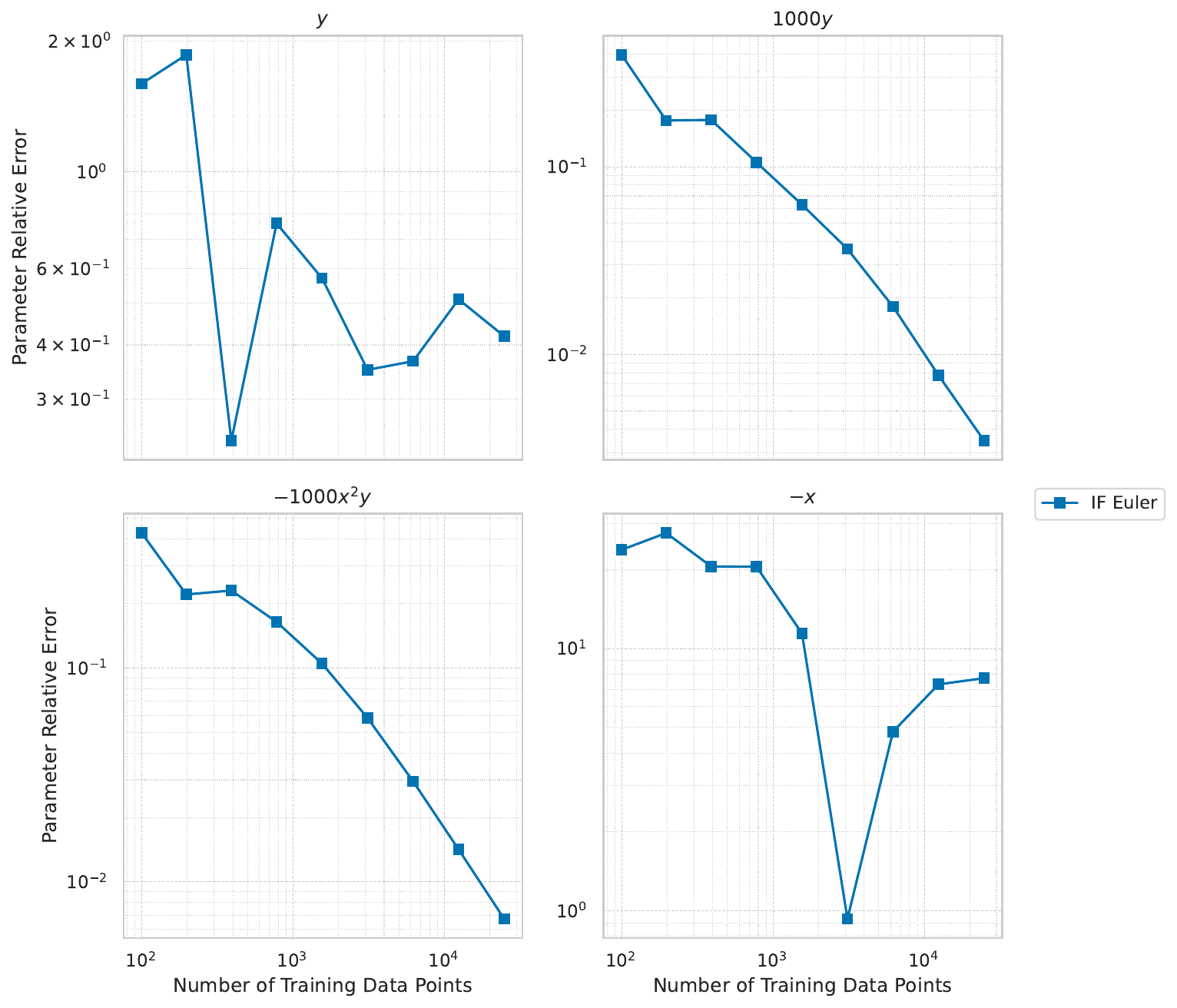}
    \caption{For the stiff Van der Pol model, we plot the fractional parameter relative error (not a percentage) against the number of training data points for the exponential integrating factor Euler integration scheme. Each subfigure corresponds to a specific term in the ODE, as indicated by the title of each subfigure.}
    \label{fig:Error_vs_n_Van_der_Pol}
\end{figure*}

Table 11 summarizes the equations recovered for different amounts of training data, and Figure \ref{fig:Error_vs_n_Van_der_Pol} shows the impact of the amount of training data on the fractional relative error of the parameters. Compared to the implicit solvers, the exponential integrating factor Euler method is remarkably inexpensive, needing only the matrix exponential calculation and no iterations. In both this case and Example 3, the performance of the method is limited by its first-order accuracy. Nonetheless, it was the only method to achieve stable training of the extremely stiff Van der Pol oscillator, maintaining stability regardless of the density of the time points in the training data.

\begin{turnpage}
\begin{table*}\centering
\ra{1.3}
\begin{tabular}{@{}rrrrcrrrcrrr@{}}\toprule
& \multicolumn{3}{c}{ 
    $\begin{aligned}
        x' &= y, \\
        y' &= 1000y - 1000x^2 y - x
    \end{aligned}$
    } \\
\cmidrule{2-4} 
& $n$ & & Equation Learned \\ \midrule
IF Euler Method\\

& $100$ & $x'$ =& $ 0.016409480199 x^{3} + 0.911941376215 x^{2} y - 0.17299439144 x^{2} - 0.0171253689169 x y^{2} - 0.0362733985181 x y  $\\
& & & $ + 0.761896059453 x - 5.24741612664 \cdot 10^{-7} y^{3} + 0.000324966428289 y^{2} - 0.596600534186 y + 1.24235801163  $\\
&    & $y'$ =& $ - 70.8004438853 x^{3} - 573.019975103 x^{2} y + 53.5327648375 x^{2} - 4.57595727557 x y^{2} + 0.227037182026 x y   $\\
& & & $ - 24.7983248216 x + 9.21699884052 \cdot 10^{-5} y^{3} - 0.00171270625655 y^{2} + 604.86402539 y + 0.719210330339  $\\

& $391$ & $x'$ =& $ 0.0211730387409 x^{3} + 1.20980396753 x^{2} y - 7.22179574621 \cdot 10^{-5} x^{2} - 0.00165302322322 x y^{2} - 0.0224633531355 x y  $\\
& & & $ + 0.163718164477 x - 7.31639160914 \cdot 10^{-7} y^{3} - 4.46875717674 \cdot 10^{-5} y^{2} + 0.759089643551 y + 0.201237070967  $\\
&    & $y'$ =& $ - 16.4449618608 x^{3} - 770.101997932 x^{2} y - 0.953912006083 x^{2} - 0.347042946264 x y^{2} + 0.148140410678 x y   $\\
& & & $ + 19.556363803 x + 0.000343170667117 y^{3} - 0.00228811438512 y^{2} + 822.181413496 y - 0.0116103130193  $\\

& $1555$ & $x'$ =& $ 0.0100174690176 x^{3} + 0.712182636703 x^{2} y - 0.000614565952893 x^{2} + 0.00133194329773 x y^{2} - 0.238901572891 x y  $\\
& & & $ + 0.00515017749321 x - 1.74495546916 \cdot 10^{-6} y^{3} - 0.000162693437889 y^{2} + 0.430535679974 y - 0.0215290827083  $\\
&    & $y'$ =& $ - 9.51940035315 x^{3} - 895.264436765 x^{2} y - 0.018305647924 x^{2} - 0.223675012257 x y^{2} - 0.831705324323 x y   $\\
& & & $ + 10.4298157344 x + 0.000139195779804 y^{3} + 0.000223913579689 y^{2} + 937.125035808 y + 0.396717128302  $\\

& $6213$ & $x'$ =& $ - 0.000116663678229 x^{3} + 0.362600279972 x^{2} y + 0.000560331851059 x^{2} + 0.00205984508829 x y^{2} + 0.152216663444 x y  $\\
& & & $ + 0.0152545009022 x - 1.84238203485 \cdot 10^{-6} y^{3} - 0.000237921205629 y^{2} + 0.633072710794 y + 0.00263343834046  $\\
&    & $y'$ =& $ - 3.18073321356 x^{3} - 970.482157645 x^{2} y - 0.057230172902 x^{2} - 0.0613306365972 x y^{2} - 0.104973101679 x y   $\\
& & & $ - 5.82202336588 x + 3.79454069868 \cdot 10^{-5} y^{3} + 5.8956384511 \cdot 10^{-5} y^{2} + 981.999703911 y + 0.348682372821  $\\

& $24849$ & $x'$ =& $ - 0.00149446009016 x^{3} + 0.185352061211 x^{2} y - 0.000178333303551 x^{2} + 0.000822839345741 x y^{2} + 0.0785771902137 x y  $\\
& & & $ + 0.0109134387138 x - 6.07013448293 \cdot 10^{-7} y^{3} - 0.000127572396163 y^{2} + 0.580986668163 y + 0.00196906484966  $\\
&    & $y'$ =& $ - 0.804037752224 x^{3} - 993.279014273 x^{2} y - 0.0226224755099 x^{2} - 0.0108890831094 x y^{2} - 0.0207393119928 x y    $\\
& & & $ - 8.70974855452 x + 5.61901765403 \cdot 10^{-6} y^{3} - 2.15618338252 \cdot 10^{-5} y^{2} + 996.524709995 y + 0.0827071351493  $\\

\bottomrule
\end{tabular}
\caption{Comparison of recovered equations for the stiff Van der Pol model using the exponential integrating factor Euler method for a selected number of training points ($n$)}
\end{table*}
\end{turnpage}

\FloatBarrier
\clearpage
\section{Conclusion}

This paper presents a new method for training stiff neural ODEs based on the explicit exponential integrating factor Euler method. Our prior study (Ref.~\onlinecite{fronk2024trainingstiffneuralordinary}) pioneered the successful training of stiff neural ODEs, showcasing robust training and accurate recovery of stiff dynamics using single-step implicit methods such as backward Euler, trapezoidal method, Radau3, and Radau5. These implicit approaches handle stiffness effectively but involve solving a nonlinear system at each time step, which can be computationally intensive if many iterations are needed for convergence. Convergence challenges intensify for these methods as the system's dimensionality grows, making their use more problematic. The explicit exponential integrating factor Euler method offers an advantage over implicit methods by needing just one matrix exponential computation per integration step, whose cost scales at \( \mathcal{O}(n^2) \) in the best case, potentially offering a more cost-effective alternative to solving nonlinear equations. 

Our results demonstrate that the IF Euler method provides stability and cost-efficiency compared to implicit methods, which struggled with stability when learning the stiff Van der Pol oscillator. In contrast, the IF Euler method successfully learned the model at large step sizes without stability issues. The findings also show that the IF Euler method is hindered by its first-order accuracy, leading to reduced precision when training neural ODEs with a fixed step size. A higher-order method would better minimize errors, though we were unable to identify a suitable higher-order explicit exponential integration method that worked (see Ref.~\onlinecite{fronk2024performanceevaluationsinglestepexplicit}). Despite its stability, its first-order accuracy results in too many steps in adaptive step size configurations, limiting its practicality for use in adaptive ODE solvers for neural ODE training. 

Our most notable finding was the exceptional performance in learning stiff linear systems of ODEs. The matrix exponential nearly perfectly solves the linear ODEs, making the exponential integration approach particularly powerful for recovering accurate linear equations. Examples 1 and 2 demonstrate that this approach scales effectively as the dimensionality of the linear system increases. The IF Euler method excels in accuracy when training data is very limited and competes with Radau5 as the dataset size grows. 

The exploration of explicit exponential integration methods provides a pathway for creating efficient, stable solvers for training stiff neural ODEs. This paves the way for applying these solvers to stiff neural PDEs, MeshGraphNets, physics-informed neural networks, and other data-driven methods that require repeatedly solving and differentiating through a differential equation. 

\section{Acknowledgements}

This work was supported in part by NSF awards CNS-1730158, ACI-1540112, ACI-1541349, OAC-1826967, OAC-2112167, CNS-2120019, the University of California Office of the President, and the University of California San Diego's California Institute for Telecommunications and Information Technology/Qualcomm Institute. Thanks to CENIC for the 100Gbps networks. The content of the information does not necessarily reflect the position or the policy of the funding agencies, and no official endorsement should be inferred.  The funders had no role in study design, data collection and analysis, decision to publish, or preparation of the manuscript.

\FloatBarrier
%\nocite{*}
\bibliography{main}% Produces the bibliography via BibTeX.

\end{document}